\input amstex
\documentstyle{amsppt}
%----------------------------------------------------------------
% Title:     Two and three descent for elliptic curves associated 
%            with perfect cuboids.
% Author:    Ruslan Sharipov
% Comments:  AmSTeX, 37 pages, amsppt style
% MSC-class: 11G05, 14H52, 11D25, 11D72  
%----------------------------------------------------------------
%           Replacement for output macro definition
%
\catcode`@=11
\redefine\output@{%
  \def\break{\penalty-\@M}\let\par\endgraf
  \ifodd\pageno\global\hoffset=105pt\else\global\hoffset=8pt\fi  
  \shipout\vbox{%
    \ifplain@
      \let\makeheadline\relax \let\makefootline\relax
    \else
      \iffirstpage@ \global\firstpage@false
        \let\rightheadline\frheadline
        \let\leftheadline\flheadline
      \else
        \ifrunheads@ %\let\makefootline\relax
        \else \let\makeheadline\relax
        \fi
      \fi
    \fi
    \makeheadline \pagebody \makefootline}%
  \advancepageno \ifnum\outputpenalty>-\@MM\else\dosupereject\fi
}
\def\Beta{\mathchar"0\hexnumber@\rmfam 42}
\catcode`\@=\active
%----------------------------------------------------------------
\nopagenumbers
\chardef\textvolna='176

\chardef\bigalpha='013
\def\negskp{\hskip -2pt}
\def\Img{\operatorname{Im}}
\def\Ker{\operatorname{Ker}}

\def\compos{\,\raise 1pt\hbox{$\sssize\circ$} \,}

\chardef\degree="5E
%\def\id{\operatorname{id}}
%\font\eightrm=cmr8
%\def\LT{\operatorname{\text{\eightrm LT}}}
%\def\LM{\operatorname{\text{\eightrm LM}}}
%\def\LC{\operatorname{\text{\eightrm LC}}}
\accentedsymbol\tildepsi{\kern 3.2pt\tilde{\kern -3.2pt\psi}}
\accentedsymbol\ssizetildepsi{\ssize\kern 2.5pt\tilde{\kern -2.5pt\psi}}
%\accentedsymbol\hatgamma{\kern 2pt\hat{\kern -2pt\gamma}}
%\accentedsymbol\checkgamma{\kern 2.5pt\check{\kern -2.5pt\gamma}}
\def\blue#1{#1}

\catcode`#=11\def\diez{#}\catcode`#=6
\catcode`&=11\catcode`&=4
\catcode`_=11\def\podcherkivanie{_}\catcode`_=8
\catcode`~=11\def\volna{~}\catcode`~=\active
\def\mycite#1{\cite{\blue{#1}}\immediate\special{ps:
     ShrHPSdict begin /ShrBORDERthickness 0 def}}

\def\mytag#1{%
    \tag#1}
\def\mythetag#1{\thetag{\blue{#1}}\immediate\special{ps:
     ShrHPSdict begin /ShrBORDERthickness 0 def}}
\def\myrefno#1{\no#1}
\def\myhref#1#2{\blue{#2}\immediate\special{ps:
     ShrHPSdict begin /ShrBORDERthickness 0 def}}
\def\myEarXivlink{\myhref{http://arXiv.org}{http:/\negskp/arXiv.org}}

\def\mytheorem#1{\csname proclaim\endcsname{Theorem #1}}
\def\mytheoremwithtitle#1#2{\csname proclaim\endcsname{Theorem #1#2}}
\def\mythetheorem#1{\blue{#1}\immediate\special{ps:
     ShrHPSdict begin /ShrBORDERthickness 0 def}}
\def\mylemma#1{\csname proclaim\endcsname{Lemma #1}}
\def\mylemmawithtitle#1#2{\csname proclaim\endcsname{Lemma #1#2}}
\def\mythelemma#1{\blue{#1}\immediate\special{ps:
     ShrHPSdict begin /ShrBORDERthickness 0 def}}
\def\mycorollary#1{\csname proclaim\endcsname{Corollary #1}}

\def\mydefinition#1{\definition{Definition #1}}

\def\myconjecture#1{\csname proclaim\endcsname{Conjecture #1}}
\def\myconjecturewithtitle#1#2{\csname proclaim\endcsname{Conjecture #1#2}}

\def\myproblem#1{\csname proclaim\endcsname{Problem #1}}
\def\myproblemwithtitle#1#2{\csname proclaim\endcsname{Problem #1#2}}

%----------------------------------------------------------------
% Cyrillic fonts definition

%\font\tencyr=wncyr10
%----------------------------------------------------------------
\pagewidth{360pt}
\pageheight{606pt}
\topmatter
\title
Two and three descent for elliptic curves associated 
with perfect cuboids.
\endtitle
\rightheadtext{Two and three descent \dots}
\author
John Ramsden, Ruslan Sharipov
\endauthor
\address Samsung Ltd\., St\. John's House, Cambridge, CB4 0ZT, UK
\endaddress
\email\myhref{mailto:jhnrmsdn\@yahoo.co.uk}{jhnrmsdn\@yahoo.co.uk}
\endemail
\address Bashkir State University, 32 Zaki Validi street, 450074 Ufa, Russia
\endaddress
\email\myhref{mailto:r-sharipov\@mail.ru}{r-sharipov\@mail.ru}
\endemail
\abstract
     A rational perfect cuboid is a rectangular parallelepiped whose edges and 
face diagonals are given by rational numbers and whose space diagonal is equal 
to unity. Recently it was shown that the Diophantine equations describing such 
a cuboid lead to a couple of parametric families of elliptic curves. Two and 
three descent methods for calculating their ranks are discussed in the present 
paper. The elliptic curves in each parametric family are subdivided into two 
subsets admitting $2$-descent and $3$-descent methods respectively.
\endabstract
\subjclassyear{2000}
\subjclass 11G05, 14H52, 11D25, 11D72\endsubjclass
\endtopmatter
%\loadbold
%\loadeufb
\TagsOnRight
\document

\head
1. Introduction.
\endhead
     Omitting the historical background of the problem, which can be found in
\mycite{1}, we proceed directly to elliptic curves in question. They are given
by the equations
$$
\xalignat 2
&\hskip -2em
2\,(w^2-1)=P_{\kern 1pt 1}\,\alpha^3,
&&2\,(w^2-1)=P_{\kern 1pt 2}\,\alpha^3
\mytag{1.1}
\endxalignat
$$
with respect to the variables $w$ and $\alpha$. The rational coefficients $P_1$
and $P_2$ in the equations \mythetag{1.1} depend on two rational parameters $b$ 
and $c$:
$$
\gathered
P_{\kern 1pt 1}=\frac{1}{2}\,(7812\,b^4\,c^4\,-216\,b^2\,c^4-52\,b^2\,c^3
+1764\,b^3\,c^4-1200\,b^4\,c^3\,-\\
-\,1848\,b^4\,c^2+720\,b^4\,c-36\,c^4\,b-1512\,b^3\,c^3-36\,c^8\,b^3
+288\,b^3\,c^2\,-\\
-\,108\,c^6\,b^2+380\,c^5\,b^2+378\,c^7\,b^3-231\,c^8\,b^4-300\,c^7\,b^4
+3906\,c^6\,b^4\,-\\
-13\,c^7\,b^2-8904\,c^5\,b^4-882\,c^6\,b^3+18\,c^6\,b-1319\,b^6\,c^8
+20952\,b^5\,c^3\,-\\
-\,11952\,b^5\,c^2+2592\,b^5\,c-48372\,b^6\,c^4+31620\,b^6\,c^3-10552\,b^6\,c^2\,+\\
+\,\,816\,b^6\,c+1494\,b^5\,c^8-5238\,b^5\,c^7-4\,c^5+7905\,b^6\,c^7
-24186\,b^6\,c^6\,+\\
+\,288\,b^6+43740\,b^6\,c^5+7686\,b^5\,c^6+576\,b^7+128\,b^8-15372\,b^5\,c^4\,-\\
-\,1080\,b^7\,c^8-3546\,b^7\,c^6+51\,c^9\,b^6+400\,b^8\,c^8-162\,c^9\,b^5
+8640\,b^7\,c^2\,-\\
-\,3456\,b^7\,c+2808\,b^7\,c^7-1560\,b^8\,c^7+3940\,b^8\,c^6+216\,c^9\,b^7
-960\,b^8\,c\,-\\
-\,6240\,b^8\,c^3+9\,c^{10}\,b^6+7880\,b^8\,c^4+4\,c^{10}\,b^8-6732\,b^8\,c^5 
+45\,c^9\,b^4\,+\\
+\,3200\,b^8\,c^2-11232\,b^7\,c^3+7092\,b^7\,c^4-18\,c^{10}\,b^7
-60\,c^9\,b^8)\,\times\\
\times\,
(b^2\,c^4-6\,b^2\,c^3+13\,b^2\,c^2-12\,b^2\,c+4\,b^2+c^2)^{-1}\!,
\endgathered
\mytag{1.2}
$$
$$
\gathered
P_{\kern 1pt 2}=\frac{b}{2}\,(832\,b^2\,c^2-1440\,b^2\,c^4-840\,b^2\,c^3
+4788\,b^3\,c^4+396\,b\,c^3\,+\\
+\,720\,b^3\,c+808\,b^4\,c^4+3032\,b^4\,c^3-2576\,b^4\,c^2
-96\,b^4\,c+448\,b^4\,-\\
-\,504\,c^4\,b-4176\,b^3\,c^3-9\,c^8\,b^3+72\,b^3\,c^2
-720\,c^6\,b^2+2288\,c^5\,b^2\,+\\
+\,1044\,c^7\,b^3-322\,c^8\,b^4+758\,c^7\,b^4+404\,c^6\,b^4
-210\,c^7\,b^2-2464\,c^5\,b^4\,-\\
-\,2394\,c^6\,b^3+72\,c^4+252\,c^6\,b+3168\,b^6\,c^8+441\,c^9\,b^5
-7056\,b^5\,c\,+\\
+\,57960\,b^6\,c^4-47232\,b^6\,c^3+25344\,b^6\,c^2
-8064\,b^6\,c-1809\,b^5\,c^8\,+\\
+\,14472\,b^5\,c^2+3951\,b^5\,c^7-72\,c^5+36\,c^6-11808\,b^6\,c^7
+1440\,b^5\,+\\
+\,28980\,b^6\,c^6-49032\,b^6\,c^5-4410\,b^5\,c^6
+8820\,b^5\,c^4-15804\,b^5\,c^3\,+\\
+\,1152\,b^6-504\,c^9\,b^6-45\,c^9\,b^3-6\,c^9\,b^4+104\,c^8\,b^2
+36\,c^{10}\,b^6\,+\\
+\,14\,c^{10}\,b^4-45\,c^{10}\,b^5-99\,c^7\,b)\,
(b^2\,c^4-6\,b^2\,c^3\,+\\
+\,13\,b^2\,c^2-12\,b^2\,c+4\,b^2+c^2)^{-1}\!.
\endgathered
\mytag{1.3}
$$
Though the equations \mythetag{1.1} arose within the same problem on perfect 
cuboids, they define two separate parametric families of elliptic curves.\par 
    The coefficients \mythetag{1.2} and \mythetag{1.3} are rational functions
of two rational parameters $b$ and $c$. Therefore, for those particular values 
of $b$ and $c$ where the common denominator of these two rational function is
nonzero, i\.\,e\. where 
$$
F=b^2\,c^4-6\,b^2\,c^3+13\,b^2\,c^2-12\,b^2\,c+4\,b^2+c^2\neq 0,
$$
their values can be represented as two irreducible fractions
$$
\xalignat 2
&\hskip -2em
P_{\kern 1pt 1}=\frac{N_1}{R_{\kern 1pt 1}},
&&P_{\kern 1pt 2}=\frac{N_2}{R_{\kern 1pt 2}}
\mytag{1.4}
\endxalignat
$$
such that $N_1\in\Bbb Z$, $N_2\in\Bbb Z$, $R_{\kern 1pt 1}\in\Bbb Z$, and 
$R_{\kern 1pt 2}\in\Bbb Z$.\par 
     Due to \mythetag{1.4} the formulas \mythetag{1.1} are transformed 
to the following ones:
$$
\xalignat 2
&\hskip -2em
2\,R_{\kern 1pt 1}\,(w^2-1)=N_1\,\alpha^3,
&&2\,R_{\kern 1pt 2}\,(w^2-1)=N_2\,\alpha^3.
\mytag{1.5}
\endxalignat
$$
Since the equations \mythetag{1.5} are very similar, we unify them by writing
the equation 
$$
\hskip -2em
2\,R\,(w^2-1)=N\,\alpha^3,
\mytag{1.6}
$$
The main goal of the present paper is to bring together some known results 
applicable to elliptic curves of the form \mythetag{1.6} thus preparing
a background for further computerized numeric search of their rational points.
Some of these points, if one is fortunate, could be responsible for perfect 
cuboids, which are wanted for centuries.\par
\head
2. Bringing to the Weierstrass form. 
\endhead
    Assuming that $R\neq 0$ and $N\neq 0$ in the equation \mythetag{1.6}, we 
substitute 
$$
\xalignat 2
&\hskip -2em
w=\frac{y}{4\,R^{\kern 0.5pt 2}\,N},
&&\alpha=\frac{x}{2\,R\,N}
\mytag{2.1}
\endxalignat
$$
into this equation. As a result we derive the equation
$$
\hskip -2em
y^2=x^3+16\,R^{\kern 0.5pt 4}\,N^{\kern 0.5pt 2}.
\mytag{2.2}
$$
The equation \mythetag{2.2} is a special case of the Weierstrass equation
$$
\hskip -2em
y^2=x^3+a\,x+b
\mytag{2.3}
$$
(see \mycite{2}) with $a=0$ and $b=16\,R^{\kern 0.5pt 4}\,N^{\kern 0.5pt 2}$. 
The transformation \mythetag{2.1} is used for to bring the cubic equation 
\mythetag{1.6} to its Weierstrass form \mythetag{2.2}. For a general cubic 
equation of two variables such a transformation bringing it to the Weierstrass 
form \mythetag{2.3} is obtained using the Nagell's algorithm (see \mycite{3}, 
\mycite{4}).\par
\head
3. The group structure. 
\endhead
     Points of an elliptic curve constitute an additive Abelian group (see \mycite{2}). 
The infinite point $P_\infty=(\infty;\,\infty)$ is usually chosen for the neutral element 
of this group: $P+P_\infty=P$. Let $P_{\kern 1pt 1}=(x_1;\,y_1)$ and 
$P_{\kern 1pt 2}=(x_2;\,y_2)$ be two points of the elliptic curve \mythetag{2.3} 
other than $P_\infty$ and let $P_{\kern 1pt 3}=(x_3;\,y_3)$ be their sum. Then, 
if $x_1\neq x_2$, the $x$-coordinate of the point $P_{\kern 1pt 3}$ is given by the 
formula 
$$
\hskip -2em
x_3=s^{\lower 1pt\hbox{$\kern 0.5pt\ssize 2$}}
-(x_1+x_2)\text{, \ where \ }s=\frac{y_1-y_2}{x_1-x_2}.
\mytag{3.1}
$$
Similarly, if $x_1\neq x_2$, the $y$-coordinate of the point $P_{\kern 1pt 3}$ 
is given by the formula 
$$
\hskip -2em
y_3=-(y_1+s\,(x_3-x_1))\text{, \ where \ }s=\frac{y_1-y_2}{x_1-x_2}
\mytag{3.2}
$$
and where $x_3$ is given by the previous formula 
\vadjust{\vskip 5pt
\hbox to 0pt{\kern 0pt \includegraphics{Cuboid_21_01.eps}
\hss}\vskip 210pt}\mythetag{3.1}.\par 
     The formulas \mythetag{3.1} and \mythetag{3.2} are applicable for 
$x_1\neq x_2$. If $x_1=x_2$, there are two options: $y_1=y_2$ and $y_1=-y_2$. 
If $x_1=x_2$ and $y_1=-y_2$, then 
$$
\hskip -2em
P_{\kern 1pt 1}+P_{\kern 1pt 2}=P_\infty
\mytag{3.3}
$$
by definition. If $x_1=x_2$ and $y_1=y_2$, the $x$-coordinate of the point 
$P_{\kern 1pt 3}$ is 
$$
\hskip -2em
x_3=s^{\lower 1pt\hbox{$\kern 0.5pt\ssize 2$}}-2\,x_1\text{, \ where \ }
s=\frac{3\,x_1^2+a}{2\,y_1}.
\mytag{3.4}
$$
Similarly, if $x_1=x_2$ and $y_1=y_2$, the $y$-coordinate of the point 
$P_{\kern 1pt 3}$ is 
$$
\hskip -2em
y_3=-(y_1+s\,(x_3-x_1))\text{, \ where \ }s=\frac{3\,x_1^2+a}{2\,y_1}
\mytag{3.5}
$$
and where $x_3$ is given by the previous \vadjust{\vskip 5pt
\hbox to 0pt{\kern 0pt \includegraphics{Cuboid_21_02.eps}
\hss}\vskip 210pt}formula \mythetag{3.4}.\par 
     The formulas \mythetag{3.1} and \mythetag{3.2} correspond to the
	case shown in Fig\.~3.1. Similarly, the formulas \mythetag{3.4} and \mythetag{3.5} 
correspond to the case shown in Fig\.~3.2. The formula \mythetag{3.3} corresponds 
to the case shown in Fig\.~3.3 above.\par
The formulas \mythetag{3.1} and \mythetag{3.2} corresponding to Fig\.~3.1 
can be transformed as
$$
\gather
\hskip -2em
x_3=\frac{x_1\,x_2^2+x_2\,x_1^2-2\,y_2\,y_1+a\,x_1+a\,x_2+2\,b}
{(x_1-x_2)^2},\qquad\quad
\mytag{3.6}\\
\vspace{1ex}
\hskip -2em
\gathered
y_3=\frac{y_2\,x_1^3-y_1\,x_2^3+4\,y_2\,b-4\,y_1\,b+y_2\,a\,x_2-y_1\,a\,x_1}
{(x_1-x_2)^3}\,+\\
+\,\frac{3\,y_2\,a\,x_1-3\,y_1\,a\,x_2+3\,y_2\,x_2\,x_1^2-3\,y_1\,x_1\,x_2^2}
{(x_1-x_2)^3}.
\endgathered
\mytag{3.7}
\endgather
$$
The formulas \mythetag{3.4} and \mythetag{3.5} corresponding to Fig\.~3.2 
can be transformed as 
$$
\pagebreak
\gather
\hskip -2em
x_3=\frac{x_1^4-2\,x_1^2\,a-8\,x_1\,b+a^2}{4\,(x_1^3+a\,x_1+b)},
\mytag{3.8}\\
\vspace{1ex}
%\displaybreak
\hskip -2em
y_3=\frac{x_1^6+5\,x_1^4\,a+20\,x_1^3\,b-5\,x_1^2\,a^2
-4\,a\,x_1\,b-a^3-8\,b^2}{8\,(x_1^3+a\,x_1+b)^2}\,y_1.
\mytag{3.9}
\endgather
$$
\par
\mydefinition{3.1} Due to the addition law given by Figs\.~3.1, 3.2, 3.3 and
by the formulas \mythetag{3.3}, \mythetag{3.6}, \mythetag{3.7}, \mythetag{3.8}, 
\mythetag{3.9} rational points of an elliptic curve $E$ given by the equation 
\mythetag{2.3}, where $a\in\Bbb Z$ and $b\in\Bbb Z$, constitute an additive 
Abelian group. This Abelian group is denoted by $E(\Bbb Q)$. 
\enddefinition
\head
4. Special points. 
\endhead
     The curve \mythetag{2.2} is a special form of the general elliptic curve
\mythetag{2.3} with $a=0$ and $b=16\,R^{\kern 0.5pt 4}\,N^{\kern 0.5pt 2}$.
Since $b=16\,R^{\kern 0.5pt 4}\,N^{\kern 0.5pt 2}$ is a square of the integer
number $4\,R^{\kern 0.5pt 2}\,N$, the curve \mythetag{2.2} has the following two 
special points:
$$
\xalignat 2
&\hskip -2em
P^{\kern 0.5pt \sssize +}_0=(0;\,4\,R^{\kern 0.5pt 2}\,N),
&&P^{\kern 0.5pt \sssize -}_0=(0;\,-4\,R^{\kern 0.5pt 2}\,N).
\mytag{4.1}
\endxalignat
$$
The special points \mythetag{4.1} of the elliptic curve \mythetag{2.2}
are shown in Fig\.~3.4. Comparing Fig\.~3.4 with Figs\.~3.2 and 3.3, we derive
$$
\xalignat 2
&\hskip -2em
2\,P^{\kern 0.5pt \sssize +}_0=P^{\kern 0.5pt \sssize +}_0
+P^{\kern 0.5pt \sssize +}_0=P^{\kern 0.5pt \sssize -}_0,
&&P^{\kern 0.5pt \sssize +}_0+P^{\kern 0.5pt \sssize -}_0
=P_\infty=0.
\mytag{4.2}
\endxalignat
$$
Using \mythetag{4.2}, one easily derives the following formulas:
$$
\xalignat 2
&\hskip -2em
3\,P^{\kern 0.5pt \sssize +}_0=0,
&&3\,P^{\kern 0.5pt \sssize -}_0=0.
\mytag{4.3}
\endxalignat
$$
The formulas \mythetag{4.3} constitute a result formulated as the following  
theorem.
\mytheorem{4.1} Let $R$ and $N$ be nonzero integers. Then the rational points
\mythetag{4.1} of the elliptic curve $E$ given by the equation \mythetag{2.2}
generate a finite subgroup of the Abelian group $E(\Bbb Q)$. This finite 
subgroup $\{P_\infty,\,P^{\kern 0.5pt \sssize +}_0,\,P^{\kern 0.5pt \sssize -}_0\}$
is isomorphic to $\Bbb Z_3$. 
\endproclaim
     Now assume that the integer number $4\,R^{\kern 0.5pt 2}\,N\neq 0$ is an 
exact cube. In this case we can write $4\,R^{\kern 0.5pt 2}\,N=8\,M^{\kern 0.5pt 3}$. 
Applying this equality to \mythetag{2.2}, we derive 
$$
\hskip -2em
y^2=x^3+64\,M^{\kern 0.5pt 6}.
\mytag{4.4}
$$
The curve \mythetag{4.4} has one more special point in addition to the points 
\mythetag{4.1}:
$$
\hskip -2em
P_{\kern 1pt 0}=(-4\,M^{\kern 0.5pt 2};\,0).
\mytag{4.5}
$$
The point \mythetag{4.5} is shown in Fig\.~3.4. Comparing Fig\.~3.4 and Fig\.~3.3,
we derive
$$
\hskip -2em
2\,P_{\kern 1pt 0}=0.
\mytag{4.6}
$$
\mytheorem{4.2} Let $M$ be a nonzero integer. Then the rational point
\mythetag{4.5} of the elliptic curve $E$ given by the equation \mythetag{4.4}
generates a finite subgroup of the Abelian group $E(\Bbb Q)$. This finite 
subgroup $\{P_\infty,\,P_{\kern 1pt 0}\}$ is isomorphic to $\Bbb Z_2$. 
\endproclaim
     The formula \mythetag{4.6} says that $P_{\kern 1pt 0}$ is a second order 
rational point of the curve \mythetag{4.4}. Similarly, the formulas \mythetag{4.3} 
say that $P^{\kern 0.5pt \sssize +}_0$ and $P^{\kern 0.5pt \sssize -}_0$ are two 
third order rational points of the curve \mythetag{2.2}. Below in sections 6 and 10
we especially study two classes of elliptic curves --- the class of curves with 
a rational point of the order $2$ and the class of curves with two rational 
points of the order $3$.
\head
4. Some classical theorems. 
\endhead
     Let $E$ be an elliptic curve given by the equation \mythetag{2.3}, where 
$a\in\Bbb Q$ and $b\in\Bbb Q$. The structure of the Abelian group $E(\Bbb Q)$ 
for such a curve is described by the following well-known theorem (see 
\mycite{5}, \mycite{6}, or \mycite{4}). 
\mytheoremwithtitle{5.1}{ (Mordell) } For an elliptic curve 
\mythetag{2.3} with $a\in\Bbb Q$ and $b\in\Bbb Q$ the group of its rational
points $E(\Bbb Q)$ is finitely generated. 
\endproclaim
     Rational points of finite order constitute a finite subgroup in $E(\Bbb Q)$.
This subgroup $E_{\sssize\text{tors}}(\Bbb Q)$ is called the torsion subgroup. 
Due to Theorem~\mythetheorem{5.1} we have
$$
\hskip -2em
E(\Bbb Q)\cong\Bbb Z^r\oplus E_{\sssize\text{tors}}(\Bbb Q)\text{, \ where \ }
r<\infty.
\mytag{5.1}
$$
The integer number $0\leqslant r<\infty$ in \mythetag{5.1} is called the 
{\it rank\/} of an elliptic curve. The structure of the torsion subgroup 
$E_{\sssize\text{tors}}(\Bbb Q)$ in the formula \mythetag{5.1} is described 
by the following theorem (see \mycite{6}).
\mytheoremwithtitle{5.2}{ (Mazur) } For an elliptic curve 
\mythetag{2.3} with $a\in\Bbb Q$ and $b\in\Bbb Q$ if the torsion subgroup of its 
rational points $E_{\sssize\text{tors}}(\Bbb Q)$ is not trivial, then it is isomorphic 
to $\Bbb Z_m$, where $m$ is one of the numbers $2,\,3,\,4,\,5,\,6,\,7,\,8,\,9,\,10,
\,12$, or it is isomorphic to $\Bbb Z_m\times\Bbb Z_2$, where $m$ is one of the numbers 
$2,\,\,4,\,6,\,\,8$.
\endproclaim
     Note that the elliptic curves associated with cuboids are not
general elliptic curves. They are given by the equation \mythetag{2.2}, which corresponds 
to $a=0$ in \mythetag{2.3}. For such curves in \mycite{6} we find the following theorem.
\mytheorem{5.3} For an elliptic curve $y^2=x^3+b$, where $b$ is a sixth-power free
integer, if the torsion subgroup of its rational points $E_{\sssize\text{tors}}(\Bbb Q)$ 
is not trivial, then it is isomorphic to $\Bbb Z_m$, where $m$ is one of the numbers 
$2,\,3,\,6$. The option $m=2$ corresponds to the case where $b$ is a cube different
from $1$. The option $m=3$ corresponds to the case where $b$ is a square different
from $1$ or $b=-2^{\kern 0.5pt 4}\,3^3$. And finally, the option $m=6$ corresponds 
to the case where $b=1$. 
\endproclaim
\noindent
     In general case the parameter $b=16\,R^{\kern 0.5pt 4}\,N^{\kern 0.5pt 2}$ in 
\mythetag{2.2} is not sixth-power free. However, one can bring it to a sixth-power 
free form using the transformation
$$
\xalignat 3
&\hskip -2em
x\to u^2\,x,
&&y\to u^3\,y,
&&b\to u^6\,b.
\quad
\mytag{5.2}
\endxalignat
$$
Theorem~\mythetheorem{5.3} along with the transformation \mythetag{5.2} means that 
the order of a point $P\in E_{\sssize\text{tors}}(\Bbb Q)$ in our case is not 
greater than $6$. The option $E_{\sssize\text{tors}}(\Bbb Q)\cong\Bbb Z_2$ is excluded 
for the curve \mythetag{2.2} due to Theorem~\mythetheorem{4.1}. Thus we have the 
following \nolinebreak result. 
\mytheorem{5.4} For an elliptic curve of the form \mythetag{2.2} associated 
with cuboids the torsion subgroup of its rational points $E_{\sssize\text{tors}}(\Bbb Q)$
is isomorphic to either $\Bbb Z_3$ or $\Bbb Z_6$. 
\endproclaim
     The case $E_{\sssize\text{tors}}(\Bbb Q)\cong\Bbb Z_6$ is the most simple. 
In this case Theorem~\mythetheorem{5.3} says that the curve \mythetag{2.2} 
is isomorphic to the curve 
$$
\hskip -2em
y^2=x^3+1. 
\mytag{5.3}
$$
Applying the transformation \mythetag{5.2} backward to the curve \mythetag{5.3}, we derive
the equality $16\,R^{\kern 0.5pt 4}\,N^{\kern 0.5pt 2}=u^6$. This equality means that the
curve \mythetag{5.3} corresponds to the case where $4\,R^{\kern 0.5pt 2}\,N\neq 0$ is an 
exact cube, i\.\,e\. $4\,R^{\kern 0.5pt 2}\,N=8\,M^{\kern 0.5pt 3}$ and $u=2\,M$. \pagebreak
This case was considered above in Theorem~\mythetheorem{4.2}.\par
     Rational points of the curve \mythetag{5.3} were studied by Euler in 1738. He has proved
the following theorem. 
\mytheoremwithtitle{5.5}{ (Euler) } If $x$ and $y$ are positive rational numbers 
satisfying the equation \mythetag{5.3}, then $x=2$ and $y=3$.
\endproclaim
     Euler's theorem~\mythetheorem{5.5} is equivalent to the following 
proposition which is a modern version of Euler's theorem~\mythetheorem{5.5}. 
\mytheorem{5.6} The rank of the curve \mythetag{5.3} is equal to zero. Its 
group of rational points $E(\Bbb Q)$ is the group of six elements
$$
E(\Bbb Q)=\{(\infty;\,\infty),\,(2;\,3),\,(0;\,1),\,(-1;\,0),\,
(0;\,-1),\,(2;\,-3)\}=E_{\sssize\text{tors}}(\Bbb Q)\cong\Bbb Z_6.
$$
\endproclaim
Applying the transformation \mythetag{5.2} backward to the curve \mythetag{5.3} and using
Theorem~\mythetheorem{5.6}, we can strengthen Theorem~\mythetheorem{4.2} in the following
way. 
\mytheorem{5.7} If the nonzero integer number $4\,R^{\kern 0.5pt 2}\,N$ is an exact cube,
i\.\,e\. if $4\,R^{\kern 0.5pt 2}\,N=8\,M^{\kern 0.5pt 3}$, then the rank of the curve
\mythetag{2.2}, which is written as \mythetag{4.4} in this case, is equal to zero. The
group of rational points $E(\Bbb Q)$ of the curve \mythetag{4.4} is finite and is composed 
by the following six elements:
$$
\xalignat 3
&P_\infty=(\infty;\,\infty), 
&&P^{\kern 0.5pt \sssize +}_*=(8\,M^{\kern 0.5pt 2};\,24\,M^{\kern 0.5pt 3}),
&&P^{\kern 0.5pt \sssize +}_0=(0;\,8\,M^{\kern 0.5pt 3}),\\
&P_{\kern 1pt 0}=(-4\,M^{\kern 0.5pt 2};\,0), 
&&P^{\kern 0.5pt \sssize -}_0=(0;\,-8\,M^{\kern 0.5pt 3}),
&&P^{\kern 0.5pt \sssize -}_*=(8\,M^{\kern 0.5pt 2};\,-24\,M^{\kern 0.5pt 3}).
\endxalignat
$$
\endproclaim
Theorem~\mythetheorem{5.7} is immediate from Theorem~\mythetheorem{5.6}. The proof of 
Theorem~\mythetheorem{5.6}, i\.\,e\. a modern proof of Euler's theorem~\mythetheorem{5.5}, 
can be found in \mycite{7}. We reproduce this proof with minor changes below in section 9 
for the sake of completeness and as an example of using the $2$-descent method.\par
\head
6. Two descent. Isogenies and descent mappings.
\endhead
     The $2$-descent method is a method of calculating the rank of an elliptic curve.
Here it is applied to curves with at least one rational point of the order $2$. Let's 
consider a general elliptic curve $E$ in Weierstrass form \mythetag{2.3} possessing a 
rational point $P_{\kern 1pt 0}$ of the order $2$. The equality $2\,P_{\kern 1pt 0}=0$ 
means $P_{\kern 1pt 0}=-P_{\kern 1pt 0}$. According to Fig\.~3.3, for a point $P=(x;\,y)$ 
on the curve \mythetag{2.3} its opposite point is $-P=(x;-y)$. Then $P_{\kern 1pt 0}
=-P_{\kern 1pt 0}$ means that the $y$-coordinate of the point $P_{\kern 1pt 0}$ is equal 
to zero. Let's denote the $x$-coordinate of the point $P_{\kern 1pt 0}$ through $c$. As
a result we get $P_{\kern 1pt 0}=(c;\,0)$. Substituting $x=c$ and $y=0$ into \mythetag{2.3}, 
we find that $b=-c^3-a\,c$. Thus, a curve possessing a point of the order $2$ is given by 
the equation
$$
\hskip -2em
y^2=x^3+a\,x-c^3-a\,c.
\mytag{6.1}
$$
The curve \mythetag{5.3} is an example of such a curve \mythetag{6.1} where 
$a=0$ and $c=-1$. The curve \mythetag{4.4} is another example with $a=0$ and 
$c=-4\,M^{\kern 0.5pt 2}$.\par
     Along with the curve \mythetag{6.1}, we consider the other curve $\tilde E$ of 
the same sort 
$$
\hskip -2em
\tilde y^2=\tilde x^3+\tilde a\,\tilde x-\tilde c^3-\tilde a\,\tilde c,
\mytag{6.2}
$$
where $\tilde a=-4\,a-15\,c^2$ and $\tilde c=-2\,c$. The curve $\tilde E$ given by
the equation \mythetag{6.2} is called the associated curve for the curve \mythetag{6.1}.
\par
     Note that the discriminant of the polynomial in the right hand side of \mythetag{6.1} 
is given by the formula $D=-(3\,c^2+4\,a)\,(a+3\,c^2)^2$, which means that the cases 
$a=-3\,c^2$ and $a=-3\,c^2/4$ in \mythetag{6.1} correspond to singular curves. Note also that 
$a=-3\,c^2$ is equivalent to $\tilde a=-3\,\tilde c^2/4$ and $a=-3\,c^2/4$ is 
equivalent to $\tilde a=-3\,\tilde c^2$. This result is formulated as a lemma.
%\mylemma{6.1} If $a\neq -3\,c^2$ and $a\neq -3\,c^2/4$, then both curves \mythetag{6.1}
%and \mythetag{6.2} are non-singular simultaneously. 
%\endproclaim
\mylemma{6.1} Both curves \mythetag{6.1} and \mythetag{6.2} are non-singular if and only if
$a\neq -3\,c^2$ and $a\neq -3\,c^2/4$. Otherwise both of them are singular. 
\endproclaim
     Now, assuming that $a\neq -3\,c^2$ and $a\neq -3\,c^2/4$ and following \mycite{7}, 
we define a mapping $\psi\!:\,E\longrightarrow\tilde E$ by means of the formulas
$$
\gather
\hskip -2em
\tilde x=\frac{x^2-x\,c+a+3\,c^2}{x-c},
\mytag{6.3}\\
\vspace{1ex}
\hskip -2em
\tilde y=\frac{y\,(x^2-2\,x\,c-a-2\,c^2)}{(x-c)^2},
\mytag{6.4}
\endgather
$$
where $x\neq c$ and $x\neq\infty$. For the special points $P_{\kern 1pt 0}$ and 
$P_\infty$ we set by definition
$$
\xalignat 2
&\hskip -2em
\psi(P_{\kern 1pt 0})=\tilde P_\infty,
&&\psi(P_\infty)=\tilde P_\infty.
\mytag{6.5}
\endxalignat
$$
Since the formulas \mythetag{6.3} and \mythetag{6.4} are not applicable to the
points $P_{\kern 1pt 0}$ and $P_\infty$, below we shall call these two points 
the exceptional points of the curve $E$. The curve $\tilde E$ has its own 
exceptional points $\tilde P_{\kern 1pt 0}$ and $\tilde P_\infty$.
\par
     Substituting \mythetag{6.3} and \mythetag{6.4} into \mythetag{6.2} and taking 
into account \mythetag{6.1}, one can prove that the formulas \mythetag{6.3}, 
\mythetag{6.4}, \mythetag{6.5} do actually define a mapping from the curve $E$ to 
its associated curve $\tilde E$.
\mylemma{6.2} The mapping $\psi\!:\,E\longrightarrow\tilde E$ defined by the formulas
\mythetag{6.3}, \mythetag{6.4}, \mythetag{6.5} induces a homomorphism of Abelian groups
$\psi\!:\,E(\Bbb Q)\longrightarrow\tilde E(\Bbb Q)$. 
\endproclaim
     Lemma~\mythelemma{6.2} implicitly assumes that the non-singularity conditions
$a\neq -3\,c^2$ and $a\neq -3\,c^2/4$ from Lemma~\mythelemma{6.1} are fulfilled.
The proof of Lemma~\mythelemma{6.2} is pure calculations using the formulas 
\mythetag{6.3}, \mythetag{6.4}, \mythetag{6.5} along with the formulas \mythetag{3.3}, 
\mythetag{3.6}, \mythetag{3.7}, \mythetag{3.8}, \mythetag{3.9}.\par
     The curve \mythetag{6.2} is of the same sort as the curve \mythetag{6.1}. Therefore
we can define a mapping $\tildepsi\!:\,\tilde E\longrightarrow E$ by means of the 
formulas 
$$
\gather
\hskip -2em
x=\frac{\tilde x^2-\tilde x\,\tilde c+\tilde a+3\,\tilde c^2}
{4\,(\tilde x-\tilde c)},
\mytag{6.6}\\
\vspace{1ex}
\hskip -2em
y=\frac{\tilde y\,(\tilde x^2-2\,\tilde x\,\tilde c-\tilde a-2\,\tilde c^2)}
{8\,(\tilde x-\tilde c)^2},
\mytag{6.7}
\endgather
$$
where $\tilde x\neq\tilde c$ and $\tilde x\neq\infty$. For the exceptional points 
$\tilde P_{\kern 1pt 0}=(\tilde c;\,0)$ and $\tilde P_\infty$ we set 
$$
\xalignat 2
&\hskip -2em
\tildepsi(\tilde P_{\kern 1pt 0})=P_\infty,
&&\tildepsi(\tilde P_\infty)=P_\infty.
\mytag{6.8}
\endxalignat
$$
Substituting \mythetag{6.6} and \mythetag{6.7} into \mythetag{6.1} and taking into 
account \mythetag{6.2}, one can verify that the formulas \mythetag{6.6}, 
\mythetag{6.7}, \mythetag{6.8} define a mapping from $\tilde E$ to $E$.
\mylemma{6.3} The mapping $\tildepsi\!:\,\tilde E\longrightarrow E$ defined by the 
formulas \mythetag{6.6}, \mythetag{6.7}, \mythetag{6.8} induces the homomorphism 
of Abelian groups $\tildepsi\!:\,\tilde E(\Bbb Q)\longrightarrow E(\Bbb Q)$. 
\endproclaim
     Lemma~\mythelemma{6.3} is similar to Lemma~\mythelemma{6.2}. Its proof
is also pure calculations.\par
     Let's consider the composite mapping $\tildepsi\compos\psi\!:\,E\longrightarrow E$.
Due to Lemmas~\mythelemma{6.2} and \mythelemma{6.3} it induces an endomorphism of
the Abelian group $E(\Bbb Q)$. 
\mylemma{6.4} The endomorphism $\tildepsi\compos\psi\!:\,E(\Bbb Q)\longrightarrow 
E(\Bbb Q)$ coincides with the doubling endomorphism, i\.\,e\. $\tildepsi\compos\psi(P)
=2\,P$ for any $P\in E(\Bbb Q)$.
\endproclaim
     The proof of Lemma~\mythelemma{6.4} is pure calculations by substituting 
\mythetag{6.3} and \mythetag{6.4} into the formulas \mythetag{6.6} and 
\mythetag{6.7}. The formulas \mythetag{6.5} and \mythetag{6.8} are also used in the
case of exceptional points.\par
      The mappings $\psi\!:\,E\longrightarrow\tilde E$ and 
$\tildepsi\!:\,\tilde E\longrightarrow E$ defined above are called $2$-isogenies
(see \mycite{8}). Due to Lemma~\mythelemma{6.4} the isogeny $\tildepsi$ is dual 
to the isogeny $\psi$ (see \mycite{9}).\par
     Let $\Bbb Q^*$ be the set of all nonzero rational numbers. This set possesses the
structure of a multiplicative Abelian group. Through $\Bbb Q^{*2}$ we denote the set
of all nonzero rational numbers which are squares. Then $\Bbb Q^{*2}$ is a subgroup
of $\Bbb Q^*$ and one can define the factor group $\Bbb Q^*\!/\Bbb Q^{*2}$. Now we 
define a mapping $\alpha\!:\,E\longrightarrow\Bbb Q^*\!/\Bbb Q^{*2}$. According
to \mycite{7}, this mapping is defined by means of the formula 
$$
\hskip -2em
\alpha(P)=x-c\text{\ \ for \ }x\neq c\text{\ \ and \ }x\neq\infty.
\mytag{6.9}
$$
The exceptional points $P_\infty$ and $P_{\kern 1pt 0}$ are treated separately.
For them we set
$$
\xalignat 2
&\hskip -2em
\alpha(P_\infty)=1,
&&\alpha(P_{\kern 1pt 0})=a+3\,c^2.
\mytag{6.10}
\endxalignat
$$
In order to relate \mythetag{6.10} with \mythetag{6.9} one should write the 
equation \mythetag{6.1} as 
$$
\align
&\hskip -2em
x-c=\frac{y^2\,(x-c)^2}{(x^2+x\,c+c^2+a)^2}\cdot
\frac{x^2+x\,c+c^2+a}{x^2-2\,x\,c+c^2},
\mytag{6.11}\\
&\hskip -2em
x-c=\frac{y^2}{(x^2+x\,c+c^2+a)^2}\cdot(x^2+x\,c+c^2+a).
\mytag{6.12}
\endalign
$$
Square factors are neglected modulo $\Bbb Q^{*2}$ in $\Bbb Q^*\!/\Bbb Q^{*2}$. 
Therefore, due to \mythetag{6.9} the formulas \mythetag{6.11} and \mythetag{6.12} 
are equivalent to the following relationships:
$$
\align
&\hskip -2em
\alpha(P)=\frac{x^2+x\,c+c^2+a}{x^2-2\,x\,c+c^2},
\mytag{6.13}\\
\vspace{1ex}
&\hskip -2em
\alpha(P)=x^2+x\,c+c^2+a.
\mytag{6.14}
\endalign
$$
Setting $x\to\infty$ in \mythetag{6.13} and $x\to c$ in \mythetag{6.14}, we 
obtain the formulas \mythetag{6.10}. 
\mylemma{6.5} The mapping $\alpha\!:\,E\longrightarrow\Bbb Q^*\!/\Bbb Q^{*2}$
defined by the formulas \mythetag{6.9} and \mythetag{6.10} induces the homomorphism 
of Abelian groups $\alpha\!:\,E(\Bbb Q)\longrightarrow\Bbb Q^*\!/\Bbb Q^{*2}$. 
\endproclaim
     In order to prove Lemma~\mythelemma{6.5} assume that $P_{\kern 1pt 1}
=(x_1;\,y_1)$ and $P_{\kern 1pt 2}=(x_2;\,y_2)$ are two non-exceptional rational 
points of the curve \mythetag{6.1} and assume that $P_{\kern 1pt 3}=(x_3;\,y_3)$
is their sum. Then \mythetag{6.9} yields $\alpha(P_{\kern 1pt 3})=x_3-c$. Applying 
the formula \mythetag{3.6}, where $b=-c^3-a\,c$ in the case of the curve 
\mythetag{6.1}, we derive
$$
\gather
\alpha(P_{\kern 1pt 3})=(x_1\,x_2^2+x_2\,x_1^2+a\,x_2-2\,c^3-2\,a\,c-2\,y_2\,y_1\,+\\
\displaybreak
+\,a\,x_1-c\,x_2^2+2\,c\,x_2\,x_1-c\,x_1^2)\,(x_1-x_2)^{-2}.
\endgather
$$
The square factor $(x_1-x_2)^{-2}$ is inessential modulo $\Bbb Q^{*2}$ in 
$\Bbb Q^*\!/\Bbb Q^{*2}$. Therefore 
$$
\hskip -2em
\gathered
\alpha(P_{\kern 1pt 3})=x_1\,x_2^2+x_2\,x_1^2+a\,x_2-2\,c^3-2\,a\,c-2\,y_2\,y_1\,+\\
+\,a\,x_1-c\,x_2^2+2\,c\,x_2\,x_1-c\,x_1^2.
\endgathered
\mytag{6.15}
$$
Now let's calculate the product $\alpha(P_{\kern 1pt 1})\,\alpha(P_{\kern 1pt 2})$. 
Applying \mythetag{6.9}, we get 
$$
\hskip -2em
\alpha(P_{\kern 1pt 1})\,\alpha(P_{\kern 1pt 2})=(x_1-c)\,(x_2-c).
\mytag{6.16}
$$
The right hand sides of \mythetag{6.15} and \mythetag{6.16} look quite 
different. However, in $\Bbb Q^*\!/\Bbb Q^{*2}$ we can multiply the right 
hand side of \mythetag{6.16} by a square factor, which is inessential:
$$
\hskip -2em
\alpha(P_{\kern 1pt 1})\,\alpha(P_{\kern 1pt 2})=(x_1-c)\,(x_2-c)
\biggl(\frac{y_1}{x_1-c}-\frac{y_2}{x_2-c}\biggr)^{\lower 2pt
\hbox{$\kern -1.5pt\ssize 2$}}.
\mytag{6.17}
$$
Expanding the right hand side of \mythetag{6.17} and transforming it with the 
use of the curve equation \mythetag{6.1}, we derive $\alpha(P_{\kern 1pt 1}+
P_{\kern 1pt 2})=\alpha(P_{\kern 1pt 3})=\alpha(P_{\kern 1pt 1})
\,\alpha(P_{\kern 1pt 2})$. Similar tricks are used for proving this equality
in the case of exceptional points $P_\infty$ and $P_{\kern 1pt 0}$.\par
     The curve \mythetag{6.2} is similar to the curve \mythetag{6.1}. Therefore
we can define a mapping $\tilde\alpha\!:\,E\longrightarrow\Bbb Q^*\!/\Bbb Q^{*2}$
using formulas similar to \mythetag{6.9} and \mythetag{6.10}:
$$
\hskip -2em
\tilde\alpha(\tilde P)=\tilde x-\tilde c\text{\ \ for \ }\tilde x\neq\tilde c
\text{\ \ and \ }\tilde x\neq\infty.
\mytag{6.18}
$$
The exceptional points $\tilde P_\infty$ and $\tilde P_{\kern 1pt 0}$ are treated 
separately. For them we set
$$
\xalignat 2
&\hskip -2em
\tilde\alpha(\tilde P_\infty)=1,
&&\tilde\alpha(\tilde P_{\kern 1pt 0})=\tilde a+3\,\tilde c^2.
\mytag{6.19}
\endxalignat
$$
\mylemma{6.6} The mapping $\tilde\alpha\!:\,\tilde E\longrightarrow
\Bbb Q^*\!/\Bbb Q^{*2}$ defined by the formulas \mythetag{6.18} and \mythetag{6.19} 
induces a homomorphism of Abelian groups 
$\tilde \alpha\!:\,\tilde E(\Bbb Q)\longrightarrow\Bbb Q^*\!/\Bbb Q^{*2}$. 
\endproclaim
     Lemma~\mythelemma{6.6} is absolutely analogous to Lemma~\mythelemma{6.5}. For
this reason we do not provide a proof of this lemma.\par
      The mappings $\alpha\!:\,E\longrightarrow\Bbb Q^*\!/\Bbb Q^{*2}$ and 
$\tilde\alpha\!:\,\tilde E\longrightarrow\Bbb Q^*\!/\Bbb Q^{*2}$ defined above are 
called descent mappings or, more specifically, $2$-descent mappings.\par
\mylemma{6.7} The kernel of the homomorphism $\alpha\!:\,E(\Bbb Q)\longrightarrow
\Bbb Q^*\!/\Bbb Q^{*2}$ coincides with the image of the homomorphism
$\tildepsi\!:\,\tilde E(\Bbb Q)\longrightarrow E(\Bbb Q)$.
\endproclaim
     Let $P=(x;\,y)$ be a non-exceptional rational point of the curve \mythetag{6.1}
such that $P\in\Ker\alpha$. Then, according to the formula \mythetag{6.9}, we have 
$$
\hskip -2em
\alpha(P)=x-c=\beta^{\kern 1pt 2}
\mytag{6.20}
$$
for some rational number $\beta\neq 0$. From \mythetag{6.20} we derive 
$x=\beta^{\kern 1pt 2}+c$. Substituting $x=\beta^{\kern 1pt 2}+c$ into
the curve equation \mythetag{6.1}, we get 
$$
\hskip -2em
\frac{y^2}{\beta^{\kern 1pt 2}}=\beta^{\kern 1pt 4}+3\,\beta^{\kern 1pt 2}\,c
+a+3\,c^2.
\mytag{6.21}
$$
We need to find a rational point $\tilde P=(\tilde x;\,\tilde y)$ of the curve
\mythetag{6.2} such that $\tildepsi(\tilde P)=P$. The formula \mythetag{6.6}
yields the following equation for the coordinate $\tilde x$: 
$$
\hskip -2em
\frac{\tilde x^2-\tilde x\,\tilde c+\tilde a+3\,\tilde c^2}
{4\,(\tilde x-\tilde c)}=\beta^{\kern 1pt 2}+c.
\mytag{6.22}
$$
Since $\tilde x-\tilde c\neq 0$ in \mythetag{6.22}, we obtain a quadratic
equation for $\tilde x$: 
$$
\hskip -2em
\tilde x^2-2\,(c-2\,\beta^{\kern 1pt 2})\,\tilde x-11\,c^2-4\,a
-8\,\beta^{\kern 1pt 2}\,c=0.
\mytag{6.23}
$$
It is easy to calculate the discriminant of the quadratic equation \mythetag{6.23}: 
$$
\hskip -2em
D=16\,(\beta^{\kern 1pt 4}+3\,\beta^{\kern 1pt 2}\,c+a+3\,c^2).
\mytag{6.24}
$$
Then we apply the standard formula for the roots of a quadratic equation. It yields
$$
\hskip -2em
\tilde x=2\,\beta^{\kern 1pt 2}+c\pm\frac{\sqrt{D}}{2}.
\mytag{6.25}
$$
Comparing \mythetag{6.24} with the formula \mythetag{6.21}, we find that the formula 
\mythetag{6.25} yields two rational solutions for the quadratic equation \mythetag{6.23}:
$$
\hskip -2em
\tilde x=2\,\beta^{\kern 1pt 2}+c\pm\frac{2\,y}{\beta}. 
\mytag{6.26}
$$ 
In order to calculate the second coordinate $\tilde y$ of the point $\tilde P$ we use
the formula \mythetag{6.7}. This formula can be written in the following form:
$$
\hskip -2em
\tilde y=\frac{8\,y\,(\tilde x+2\,c)^2}{\tilde x^2+4\,\tilde x\,c+7\,c^2+4\,a}.
\mytag{6.27}
$$
Substituting \mythetag{6.26} into \mythetag{6.27} and applying the curve equation
\mythetag{6.21}, we get 
$$
\hskip -2em
\tilde y=4\,y\pm(4\,\beta^{\kern 1pt 3}+6\,c\,\beta). 
\mytag{6.28}
$$
The formulas \mythetag{6.26} and \mythetag{6.28} yield explicit expressions 
for the coordinates of the required rational point $\tilde P$ such that 
$\tildepsi(\tilde P)=P$. Thus, for a non-exceptional rational point $P$ of the 
curve $E$ we have proved that $P\in\Ker\alpha$ implies $P\in\Img\tildepsi$.
\par
     Let's proceed to exceptional points $P_\infty$ and $P_{\kern 1pt 0}$. The case
of the point $P_\infty$ is trivial since in this case $P_\infty\in\Ker\alpha$ and 
$P_\infty=\tildepsi(\tilde P_\infty)$.\par
     Assume that $P_{\kern 1pt 0}\in\Ker\alpha$. Then from \mythetag{6.10} we 
derive \ $a+3\,c^2=\beta^{\kern 1pt 2}$, where $\beta\neq 0$ is some rational 
number. Let's resolve the equality $a+3\,c^2=\beta^{\kern 1pt 2}$ with respect 
to $a$:
$$
\hskip -2em
a=\beta^{\kern 1pt 2}-3\,c^2.
\mytag{6.29}
$$
Taking into account the relationships $\tilde a=-4\,a-15\,c^2$, $\tilde c=-2\,c$
and applying the formula \mythetag{6.29} to \mythetag{6.2}, we find that the curve 
equation \mythetag{6.2} turns to
$$
\hskip -2em
\tilde y^2=(x+2\,c)(x-c-2\,\beta)(x-c+2\,\beta).
\mytag{6.30}
$$
The formula \mythetag{6.30} means that under the assumption $P_{\kern 1pt 0}\in
\Ker\alpha$ \pagebreak the curve $\tilde E$ has not only the exceptional point
$\tilde P_{\kern 1pt 0}=(-2\,c;\,0)$ but two other similar points
$$
\xalignat 2
&\hskip -2em
\tilde P_{\kern 1pt 01}=(c+2\,\beta;\,0),
&&\tilde P_{\kern 1pt 02}=(c-2\,\beta;\,0).
\mytag{6.31}
\endxalignat
$$
The non-exceptional points $\tilde P_{\kern 1pt 01}$ and $\tilde P_{\kern 1pt 02}$ 
neither can coincide with each other nor with the point $\tilde P_{\kern 1pt 0}$
since otherwise the curve $\tilde E$ would be singular 
(see Lemma~\mythelemma{6.1}). Applying \mythetag{6.6} and \mythetag{6.7} to the 
coordinates of the points \mythetag{6.31}, we derive 
$$
\xalignat 2
&\hskip -2em
\tildepsi(\tilde P_{\kern 1pt 01})=P_{\kern 1pt 0},
&&\tildepsi(\tilde P_{\kern 1pt 02})=P_{\kern 1pt 0}.
\mytag{6.32}
\endxalignat
$$
Each of the two formulas \mythetag{6.32} is sufficient to conclude that if 
$P_{\kern 1pt 0}\in\Ker\alpha$, then $P_{\kern 1pt 0}\in\Img\tildepsi$. Summarizing
the above considerations for exceptional and non-exceptional points, we derive
$\Ker\alpha\subseteq\Img\tildepsi$.\par 
     Now, conversely, assume that $P$ is a non-exceptional rational point 
such that $P\in\Img\tildepsi$. Then its coordinates $x$ and $y$ are given by the 
formulas \mythetag{6.6} and \mythetag{6.7}, where $\tilde x$ and $\tilde y$
are the coordinates of some rational point $\tilde P$ of the curve \mythetag{6.2}. 
Note that the formula \mythetag{6.6} can be written as follows: 
$$
\hskip -2em
x-c=\frac{\tilde x^2-2\,\tilde x\,c-11\,c^2-4\,a}{4\,(\tilde x+2\,c)}.
\mytag{6.33}
$$
Multiplying the numerator and the denominator of the fraction \mythetag{6.33} by
$\tilde x+2\,c$, we derive the following formula for $x-c$\,:
$$
\hskip -2em
x-c=\frac{\tilde x^3-(15\,c^2+4\,a)\,\tilde x-22\,c^3-8\,a\,c}{4\,(\tilde x+2\,c)^2}.
\mytag{6.34}
$$
Comparing the numerator of \mythetag{6.34} with the curve equation \mythetag{6.2},
where $\tilde c=-2\,c$ and $\tilde a=-4\,a-15\,c^2$, and taking into account 
\mythetag{6.9}, we transform \mythetag{6.34} as follows:
$$
\hskip -2em
\alpha(P)=x-c=\frac{\tilde y^2}{4\,(\tilde x+2\,c)^2}.
\mytag{6.35}
$$
The right hand side of \mythetag{6.35} is a square of a rational number. For this 
reason $\alpha(P)\in\Ker\alpha$. Thus, for a non-exceptional rational point $P$ of 
the curve $E$ we have proved that $P\in\Img\tildepsi$ implies $P\in\Ker\alpha$.
\par
     Let's proceed to exceptional points $P_\infty$ and $P_{\kern 1pt 0}$. 
The case of the point $P_\infty$ is trivial since in this case 
$P_\infty=\tildepsi(\tilde P_\infty)$ and $P_\infty\in\Ker\alpha$.\par
     Assume that $P_{\kern 1pt 0}\in\Img\tildepsi$. Then the coordinates 
$x=c$ and $y=0$ of the point $P_{\kern 1pt 0}$ are given by the formulas 
\mythetag{6.6} and \mythetag{6.7}, where $\tilde x$ and $\tilde y$
are the coordinates of some non-exceptional rational point $\tilde P$ of the 
curve \mythetag{6.2}. Therefore the equality $x=c$ leads to the following 
equations with respect to $\tilde x$:
$$
\hskip -2em
\frac{\tilde x^2-\tilde x\,\tilde c+\tilde a+3\,\tilde c^2}
{4\,(\tilde x-\tilde c)}=c,
\mytag{6.36}
$$
The denominator of the fraction in \mythetag{6.36} is nonzero. 
Therefore the equation \mythetag{6.36} is equivalent to a quadratic equation for 
$\tilde x$. \pagebreak Due to $\tilde a=-4\,a-15\,c^2$ and $\tilde c=-2\,c$ this 
quadratic equation can be written as   
$$
\hskip -2em
\tilde x^2-2\,\tilde x\,c-11\,c^2-4\,a=0.
\mytag{6.37}
$$
One can easily calculate the discriminant $D=16\,(a+3\,c^2)$ 
of the quadratic equation \mythetag{6.37} and derive the following explicit 
formula for its solution $\tilde x$:
$$
\hskip -2em
\tilde x=c\pm2\,\sqrt{a+3\,c^2\,}.
\mytag{6.38}
$$
The formula \mythetag{6.38} can be transformed to the following one:
$$
\hskip -2em
a+3\,c^2=\Bigl(\frac{\tilde x-c}{2}\kern 0.3pt\Bigr)^2\!.
\mytag{6.39}
$$
Comparing \mythetag{6.39} with the formula \mythetag{6.10} for 
$\alpha(P_{\kern 1pt 0})$, we derive 
$$
\hskip -2em
\alpha(P_{\kern 1pt 0})=\Bigl(\frac{\tilde x-c}{2}\kern 0.3pt\Bigr)^2\!.
\mytag{6.40}
$$
Square factors are neglected modulo $\Bbb Q^{*2}$ in $\Bbb Q^*\!/\Bbb Q^{*2}$. 
Hence the equality \mythetag{6.40} is equivalent to $\alpha(P_{\kern 1pt 0})=1$, 
which means $P_{\kern 1pt 0}\in\Ker\alpha$. We have proved that 
$P_{\kern 1pt 0}\in\Img\tildepsi$ implies $P_{\kern 1pt 0}\in\Ker\alpha$. 
Summarizing the above considerations for exceptional and non-exceptional points, 
we derive $\Img\tildepsi\subseteq\Ker\alpha$. Then, combining 
$\Img\tildepsi\subseteq\Ker\alpha$ with the previously derived inclusion 
$\Ker\alpha\subseteq\Img\tildepsi$, we conclude that $\Ker\alpha=\Img\tildepsi$, 
which completes the proof of Lemma~\mythelemma{6.7}. 
\mylemma{6.8} The kernel of the homomorphism $\tilde\alpha\!:\,\tilde E(\Bbb Q)
\longrightarrow \Bbb Q^*\!/\Bbb Q^{*2}$ coincides with the image of the 
homomorphism $\psi\!:\,E(\Bbb Q)\longrightarrow\tilde E(\Bbb Q)$.
\endproclaim
     Lemma~\mythelemma{6.8} is proved in a way similar to the above proof
of Lemma~\mythelemma{6.7}.\par
\head
7. Factor groups and the rank formula. 
\endhead
     Let's proceed to further studying the homomorphisms $\psi\!:\,E(\Bbb Q)
\longrightarrow\tilde E(\Bbb Q)$ and $\tildepsi\!:\,\tilde E(\Bbb Q)\longrightarrow 
E(\Bbb Q)$. The visual image in Fig\.~7.1 will help the reader. \vadjust{\vskip 5pt
\hbox to 0pt{\kern 0pt \includegraphics{Cuboid_21_03.eps}
\hss}\vskip 170pt}In the left hand side of this image we have two independent
inclusions $\psi(E(\Bbb Q))\subset\tilde E(\Bbb Q)$ and $\Ker\tildepsi\subset
\tilde E(\Bbb Q)$. \pagebreak In the right hand side of the image we have the series 
of inclusions $\tildepsi\compos\psi(E(\Bbb Q))\subset\tildepsi(\tilde E(\Bbb Q))
\subset E(\Bbb Q)$. From this series of inclusions we derive 
$$
\hskip -2em
E/\tildepsi(\tilde E(\Bbb Q))\cong(E/\tildepsi\compos\psi(E(\Bbb Q)))\,/
\,(\tildepsi(\tilde E(\Bbb Q))/\tildepsi\compos\psi(E(\Bbb Q)))
\mytag{7.1}
$$
(see basics of the group theory in \mycite{10}). The isomorphism \mythetag{7.1}
yields the following relationship for the indices of subgroups in \mythetag{7.1}:
$$
\hskip -2em
[E(\Bbb Q):2\,E(\Bbb Q)]=[E(\Bbb Q):\tildepsi(\tilde E(\Bbb Q))]
\cdot[\tildepsi(\tilde E(\Bbb Q)):\tildepsi\compos\psi(E(\Bbb Q))].
\mytag{7.2}
$$
In deriving the above equality \mythetag{7.2} we used the equality $2\,E(\Bbb Q)
=\tildepsi\compos\psi(E(\Bbb Q))$ provided by Lemma~\mythelemma{6.4}.\par 
     Note that the mapping $\tildepsi\!:\,\tilde E(\Bbb Q)\longrightarrow 
E(\Bbb Q)$ in Fig\.~7.1 can be treated as the surjective mapping $\tildepsi\!:\,
\tilde E(\Bbb Q)\longrightarrow\tildepsi(\tilde E(\Bbb Q))$. This surjective mapping 
can be combined with the factorization mapping $\tau\!:\,\tildepsi(\tilde E(\Bbb Q))
\longrightarrow\tildepsi(\tilde E(\Bbb Q))/\tildepsi\compos\psi(E(\Bbb Q))$:
$$
\hskip -2em
\CD
\tilde E(\Bbb Q)@>\ \ssizetildepsi\ >>\tildepsi(\tilde E(\Bbb Q))@>\ \tau\ >>
\tildepsi(\tilde E(\Bbb Q))/\tildepsi\compos\psi(E(\Bbb Q)).
\endCD
\mytag{7.3}
$$
Both mappings in \mythetag{7.3} are surjective. Therefore the composite mapping 
$\tau\compos\tildepsi$ in \mythetag{7.3} is also surjective. Its kernel is easily
calculated:
$$
\hskip -2em
\Ker(\tau\compos\tildepsi)=\psi(E(\Bbb Q))+\Ker\tildepsi.
\mytag{7.4}
$$
From \mythetag{7.3} and \mythetag{7.4} we derive the isomorphism
$$
\hskip -2em
\tilde E(\Bbb Q)/(\psi(E(\Bbb Q))+\Ker\tildepsi)\cong
\tildepsi(\tilde E(\Bbb Q))/\tildepsi\compos\psi(E(\Bbb Q))
\mytag{7.5}
$$
Along with the isomorphism \mythetag{7.5} we have the following two isomorphisms:
$$
\gather
\hskip -2em
\gathered
\tilde E(\Bbb Q)/(\psi(E(\Bbb Q))+\Ker\tildepsi)\,\cong\\
\cong\,(\tilde E(\Bbb Q)/\psi(E(\Bbb Q)))\,/\,((\psi(E(\Bbb Q))+\Ker\tildepsi)
/\psi(E(\Bbb Q))),
\endgathered
\mytag{7.6}\\
\vspace{1.5ex}
\hskip -2em
(\psi(E(\Bbb Q))+\Ker\tildepsi)/\psi(E(\Bbb Q))\cong
\Ker\tildepsi/(\Ker\tildepsi\cap\psi(E(\Bbb Q)).
\mytag{7.7}
\endgather
$$
The isomorphisms \mythetag{7.6} and \mythetag{7.7} are basic facts from the
group theory (see \mycite{10}). From \mythetag{7.5}, \mythetag{7.6}, and
\mythetag{7.7} we derive the equalities
$$
\gather
\hskip -2em
[\tildepsi(\tilde E(\Bbb Q)):\tildepsi\compos\psi(E(\Bbb Q))]
=[\tilde E(\Bbb Q):(\psi(E(\Bbb Q))+\Ker\tildepsi)],
\mytag{7.8}\\
\vspace{1ex}
\hskip -2em
[\tilde E(\Bbb Q):(\psi(E(\Bbb Q))+\Ker\tildepsi)]=
\frac{[\tilde E(\Bbb Q):\psi(E(\Bbb Q))]}{[(\psi(E(\Bbb Q))+\Ker\tildepsi)
:\psi(E(\Bbb Q))]},
\mytag{7.9}\\
\vspace{2ex}
\hskip -2em
[(\psi(E(\Bbb Q))+\Ker\tildepsi):\psi(E(\Bbb Q))]
=[\Ker\tildepsi:(\Ker\tildepsi\cap\psi(E(\Bbb Q))].
\mytag{7.10}
\endgather
$$
Now, combining \mythetag{7.2} with \mythetag{7.8}, \mythetag{7.9}, and
\mythetag{7.10}, we obtain the equality
$$
\hskip -2em
[E(\Bbb Q):2\,E(\Bbb Q)]=\frac{[E(\Bbb Q):\tildepsi(\tilde E(\Bbb Q))]
\cdot[\tilde E(\Bbb Q):\psi(E(\Bbb Q))]}
{[\Ker\tildepsi:(\Ker\tildepsi\cap\psi(E(\Bbb Q))]}.
\mytag{7.11}
$$
\mylemma{7.1} The index $[\Ker\tildepsi:(\Ker\tildepsi\cap\psi(E(\Bbb Q))]$ equals
$1$ if and only if the curve \mythetag{6.1} has a rational point $P=(x,\,0)$ of the 
order $2$ different from its exceptional point $P_{\kern 1pt 0}=(c;\,0)$. Otherwise 
$[\Ker\tildepsi:(\Ker\tildepsi\cap\psi(E(\Bbb Q))]=2$.
\endproclaim
     Looking at \mythetag{6.8}, we see that $\Ker\tildepsi$ consists of two elements
$\tilde P_\infty$ and $\tilde P_{\kern 1pt 0}$. According to \mythetag{6.5} the point 
$\tilde P_\infty$ belongs to $\psi(E(\Bbb Q))$. The point $\tilde P_{\kern 1pt 0}=
(\tilde c;\,0)$ belongs to $\psi(E(\Bbb Q))$ if and only if there is a non-exceptional 
rational point $P=(x,y)$ of the curve \mythetag{6.1} such that $\psi(P)=\tilde P_{\kern 
1pt 0}$. Applying \mythetag{6.4} and \mythetag{6.5}, we derive 
$$
\gather
\hskip -2em
\frac{x^2-x\,c+a+3\,c^2}{x-c}=\tilde c,
\mytag{7.12}\\
\vspace{1ex}
\hskip -2em
\frac{y\,(x^2-2\,x\,c-a-2\,c^2)}{(x-c)^2}=0,
\mytag{7.13}
\endgather
$$
where $x\neq c$. Since $\tilde c=-2\,c$, the equality \mythetag{7.12} reduces to 
$$
\hskip -2em
x^2+x\,c+c^2+a=0.
\mytag{7.14}
$$ 
The polynomial $x^2+x\,c+c^2+a$ cannot vanish simultaneously
with the polynomial $x^2-2\,x\,c-a-2\,c^2$ in the numerator of the fraction in
\mythetag{7.13}. Otherwise we would have the following system of two quadratic
equations: 
$$
\left\{\aligned
&x^2+x\,c+c^2+a=0,\\
&x^2-2\,x\,c-a-2\,c^2=0.
\endaligned\right.
\mytag{7.15}
$$
Let's multiply the first equation \mythetag{7.15} by $9\,c^2-3\,x\,c+2\,a$ and 
multiply the second equation \mythetag{7.15} by $3\,x\,c-2\,a$. Then, adding the
resulting two equations, we derive 
$$
\hskip -2em
9\,c^4+15\,c^2\,a+4\,a^2=0.
\mytag{7.16}
$$
The left hand side of the equality \mythetag{7.16} factors as follows:
$$
\hskip -2em
(a+3\,c^2)\,(4\,a+3\,c^2)=0.
\mytag{7.17}
$$
Applying Lemma~\mythelemma{6.1}, we see that the equality \mythetag{7.17}
contradicts the non-singularity condition of the curves \mythetag{6.1} and
\mythetag{6.2}.\par
     The contradiction obtained means that the equation \mythetag{7.14} should
be complemented with the inequality $x^2-2\,x\,c-a-2\,c^2\neq 0$. Then from
\mythetag{7.13} we derive $y=0$, which means that $P=(x;\,y)$ is a rational
point of the order $2$. Since $x\neq c$ in \mythetag{7.12} and \mythetag{7.13},
this point does not coincide with $P_{\kern 1pt 0}$. Lemma~\mythelemma{7.1} 
is proved.\par
     Let's consider the subgroup indices $[E(\Bbb Q):\tildepsi(\tilde E(\Bbb Q))]$ 
and $[\tilde E(\Bbb Q):\psi(E(\Bbb Q))]$ in the numerator of the fraction in the
right hand side of \mythetag{7.11}. Applying Lemmas~\mythelemma{6.7} and 
\mythelemma{6.8}, for these indices we derive 
$$
\xalignat 2
&[E(\Bbb Q):\tildepsi(\tilde E(\Bbb Q))]=|\alpha(E(\Bbb Q))|,
&&[\tilde E(\Bbb Q):\psi(E(\Bbb Q))]=|\tilde\alpha(\tilde E(\Bbb Q))|.
\qquad\quad
\mytag{7.18}
\endxalignat 
$$
The number of elements in a group or, more generally, in a set $G$ is often denoted
trough $\#G$. But we prefer to use the notation $|G|$ in \mythetag{7.18} \pagebreak 
instead of $\#G$. Due to the formulas \mythetag{7.18} the formula \mythetag{7.11} 
is written in the following way:
$$
\hskip -2em
[E(\Bbb Q):2\,E(\Bbb Q)]=\frac{|\alpha(E(\Bbb Q))|
\cdot|\tilde\alpha(\tilde E(\Bbb Q))|}
{[\Ker\tildepsi:(\Ker\tildepsi\cap\psi(E(\Bbb Q))]}.
\mytag{7.19}
$$
Let's recall Mordell's theorem~\mythetheorem{5.1} and the formula \mythetag{5.1}.
From \mythetag{5.1} we derive 
$$
\hskip -2em
E(\Bbb Q)/2\,E(\Bbb Q)\cong (Z_2)^r\oplus E_2(\Bbb Q),
\mytag{7.20}
$$
where $E_2(\Bbb Q)$ is the subgroup of $E_{\sssize\text{tors}}(\Bbb Q)$ generated by
all elements of the order $2$. Applying \mythetag{7.20} to \mythetag{7.19}, we transform
the formula \mythetag{7.19} as follows:
$$
\hskip -2em
2^{\kern 1pt r}\!\cdot|E_2(\Bbb Q)|=\frac{|\alpha(E(\Bbb Q))|
\cdot|\tilde\alpha(\tilde E(\Bbb Q))|}
{[\Ker\tildepsi:(\Ker\tildepsi\cap\psi(E(\Bbb Q))]}.
\mytag{7.21}
$$
The formula \mythetag{7.21} yields a way for calculating the rank $r$ of an 
elliptic curve.\par
\head
8. Images of the descent mappings. 
\endhead
     The next step is to calculate the groups $\alpha(E(\Bbb Q))$ and 
$\tilde\alpha(\tilde E(\Bbb Q))$. They are images of the descent mappings
$\alpha\!:\,E(\Bbb Q)\longrightarrow\Bbb Q^*\!/\Bbb Q^{*2}$ and 
$\tilde \alpha\!:\,\tilde E(\Bbb Q)\longrightarrow\Bbb Q^*\!/\Bbb Q^{*2}$. 
For calculating them let's present the equation \mythetag{6.1} in the 
following form:
$$
\hskip -2em
y^2=(x-c)\,((x-c)^2+3\,c\,(x-c)+B),
\mytag{8.1}
$$
where $B=3\,c^2+a$. The integer number $B$ in \mythetag{8.1} is nonzero due to
Lemma~\mythelemma{6.1}.\par
\mylemma{8.1} Let $P$ be a rational point of the curve \mythetag{6.1}. Then 
$\alpha(P)$ is presented by some integer number being a divisor of\/ $B=3\,c^2+a$.
\endproclaim
     For the exceptional points $P_\infty$ and $P_{\kern 1pt 0}$ the
proposition of Lemma~\mythelemma{8.1} follows from the formulas \mythetag{6.10}.
Let $P=(x;\,y)$ be a non-exceptional rational point of the curve \mythetag{6.1}. 
Then $x$ and $y$ are rational numbers satisfying the equation \mythetag{8.1} and 
such that $x\neq c$. Therefore they are presented by the formulas 
$$
\xalignat 2
&\hskip -2em
x=c+\frac{m}{q_1},
&y=\frac{n}{q_2},
\mytag{8.2}
\endxalignat
$$
where $m$, $q_1$, and $q_2$ are nonzero integer numbers such that the fractions
\mythetag{8.2} are irreducible. Assume that $n\neq 0$. Upon substituting \mythetag{8.2} 
into \mythetag{8.1}, we obtain 
$$
\hskip -2em
\frac{n^2}{q_2^2}=\frac{m}{q_1}\cdot\frac{m^2+3\,c\,m\,q_1+B\,q_1^2}{q_1^2}. 
\mytag{8.3}
$$
It is easy to see that all of the three fractions in \mythetag{8.3} are irreducible.
This yields 
$$
\hskip -2em
q_2^2=q_1^3. 
\mytag{8.4}
$$
Due to \mythetag{8.4} there is an integer number $q$ such that 
$$
\xalignat 2
&\hskip -2em
q_1=q^2,
&q_2=q^3
\mytag{8.5}
\endxalignat
$$
(compare with Lemma~2.1 in \mycite{1}). Substituting \mythetag{8.5} into 
\mythetag{8.2}, we get
$$
\xalignat 2
&\hskip -2em
x=c+\frac{m}{q^2},
&y=\frac{n}{q^3}.
\mytag{8.6}
\endxalignat
$$
Similarly, substituting \mythetag{8.5} into the equation \mythetag{8.3}, we 
obtain the equation
$$
\hskip -2em
n^2=m\,(m^2+3\,c\,m\,q^2+B\,q^4). 
\mytag{8.7}
$$
The left hand side of \mythetag{8.7} is positive. Therefore we can write 
\mythetag{8.7} as 
$$
\hskip -2em
n^2=|m|\cdot|m^2+3\,c\,m\,q^2+B\,q^4|. 
\mytag{8.8}
$$
Let's denote $\beta=\gcd(|m|,|m^2+3\,c\,m\,q^2+B\,q^4|)$. The fractions in 
\mythetag{8.6} are irreducible. Therefore we have $\beta=\gcd(|m|,B)$, i\.\,e\.
$\beta$ is a divisor of $B$. Then
$$
\xalignat 2
&\hskip -2em
m=\pm\,\beta\,\mu,
&&m^2+3\,c\,m\,q^2+B\,q^4=\pm\,\beta\,\nu,
\mytag{8.9}
\endxalignat
$$ 
where $\mu>0$, $\nu>0$, and $\gcd(\mu,\nu)=1$. Applying \mythetag{8.9} to 
\mythetag{8.8}, we get 
$$
\hskip -2em
n^2=\beta^{\kern 1pt 2}\kern 0.5pt\mu\,\nu. 
\mytag{8.10}
$$
The positive integers $\mu$ and $\nu$ in \mythetag{8.10} are coprime, while their 
product is a square. Therefore both of these integer numbers are squares:
$$
\xalignat 2
&\hskip -2em
\mu=M^{\kern 1pt 2},
&&\nu=N^{\kern 1pt 2}.
\mytag{8.11}
\endxalignat
$$ 
Applying \mythetag{8.11} to \mythetag{8.9} and then applying \mythetag{8.9} to
\mythetag{8.6}, we derive
$$
\hskip -2em
x-c=\frac{m}{q^2}=\pm\frac{\beta\,\mu}{q^2}
=\pm\frac{\beta\,M^{\kern 1pt 2}}{q^2}\equiv\pm\,\beta\!\!\!\mod\Bbb Q^{*2},
\mytag{8.12}
$$
where $\beta$ is a divisor of $B$. The case $n=0$ in \mythetag{8.3} is special. In this
case 
$$
\hskip -2em
x-c=\frac{m}{q_1}=\beta
\mytag{8.13}
$$
and \mythetag{8.3} yields that $\beta$ is a rational root of the monic polynomial 
$\beta^{\kern 1pt 2}-3\,c\,\beta+B=0$ with integer coefficients. In this case the 
rational root theorem (see \mycite{11}) says that $q_1=1$, while $\beta=m$ is a 
divisor of $B$. Comparing \mythetag{8.12} and \mythetag{8.13} with \mythetag{6.9}, 
we find that Lemma~\mythelemma{8.1} is proved. 
\mylemma{8.2} Let $\tilde P$ be a point of the curve \mythetag{6.2}. Then 
$\tilde\alpha(\tilde P)$ is presented by some integer number being a divisor of\/ 
$\tilde B=3\,\tilde c^2+\tilde a$.
\endproclaim
     Lemma~\mythelemma{8.2} is similar to the previous lemma~\mythelemma{8.1}. 
Its proof is similar to the above proof of Lemma~\mythelemma{8.1}.\par
\head
9. Proof of Euler's theorem. 
\endhead
\mylemma{9.1} For the curve \mythetag{5.3} the group $\alpha(E(\Bbb Q))$ has two
elements presented by the numbers \pagebreak $1$ and $3$, i\.\,e\. $|\alpha(E(\Bbb Q))|=2$.
\endproclaim
     The curve \mythetag{5.3} is a particular case of the curve \mythetag{6.1} where 
$a=0$ and where $c=-1$. Then $B=3\,c^2+a=3$. The number $3$ has two positive 
divisors $1$ and $3$. Therefore, according to Lemma~\mythelemma{8.1}, the group 
$\alpha(E(\Bbb Q))$ has at most four elements presented by the numbers $1$, $-1$,
$3$, $-3$. The formulas \mythetag{6.10} yield
$$
\xalignat 2
&\hskip -2em
\alpha(P_\infty)=1,
&&\alpha(P_{\kern 1pt 0})=3.
\mytag{9.1}
\endxalignat
$$
The relationships \mythetag{9.1} mean that $\alpha(E(\Bbb Q))$ does actually
comprise two elements presented by the numbers $1$ and $3$. Since $-3=(-1)\cdot 3$,
it is sufficient to prove that there is no element presented by the number $-1$
in $\alpha(E(\Bbb Q))$. Assume to the contrary that $\alpha(P)=-1$ for some rational
point $P$. Due to \mythetag{9.1} we have $P\neq P_\infty$ and $P\neq P_{\kern 1pt 0}$,
i\.\,e\. $P=(x;\,y)$ is a non-exceptional point. Then we can apply \mythetag{8.2} to
its coordinates and derive \mythetag{8.3}. The curve \mythetag{5.3} has no rational 
points with $y=0$ other than $P_{\kern 1pt 0}=(-1;\,0)$. Therefore $n\neq 0$ in 
\mythetag{8.2}. Then, repeating the logic of the proof of Lemma~\mythelemma{8.1},
we come from \mythetag{8.3} to \mythetag{8.12}. But now
$$
\hskip -2em
\alpha(P)=x-c\equiv -1\!\!\!\mod\Bbb Q^{*2}.
\mytag{9.2}
$$
The formula \mythetag{9.2} means that $\beta=K^{\kern 1pt 2}$ for some nonzero
integer $K$ and that we should choose the negative sign in \mythetag{8.12}. The sign
choices in \mythetag{8.12} and \mythetag{8.9} do coincide. Therefore, from
\mythetag{8.9}, \mythetag{8.10}, and \mythetag{8.11}  we derive 
$$
\xalignat 2
&\hskip -2em
m=-M^{\kern 1pt 2}\,K^{\kern 1pt 2},
&&n^2=M^{\kern 1pt 2}\,N^{\kern 1pt 2}\,K^{\kern 1pt 4}.
\mytag{9.3}
\endxalignat
$$
Substituting \mythetag{9.3} into \mythetag{8.7}, and recalling that now $c=-1$ and 
$B=3$, we get  
$$
\hskip -2em
N^{\kern 1pt 2}\,K^{\kern 1pt 2}=
-(M^{\kern 1pt 4}\,K^{\kern 1pt 4}
+3\,M^{\kern 1pt 2}\,K^{\kern 1pt 2}\,q^2
+3\,q^4).
\mytag{9.4}
$$ 
The equality \mythetag{9.4} means that $K^{\kern 1pt 2}$ divides $3\,q^4$, i\.\,e\.
if $K^{\kern 1pt 2}\neq 1$, then $K$ and $q$ would have at least one common prime 
factor $p$, which is either $p=3$ or $p\neq 3$. Due to \mythetag{9.3} this prime 
factor $p$ of $q$ would be a prime factor of $m$ and $n$ as well. But the fractions
in \mythetag{8.1} and \mythetag{8.6} are irreducible. Therefore $K^{\kern 1pt 2}=1$
and $\beta=1$. Substituting $K^{\kern 1pt 2}=1$ into the equality \mythetag{9.4}, 
we reduce this equality to 
$$
\hskip -2em
N^{\kern 1pt 2}=-M^{\kern 1pt 4}-3\,M^{\kern 1pt 2}\,q^2-3\,q^4.
\mytag{9.5}
$$ 
The equality \mythetag{9.5} means that $N^{\kern 1pt 2}\equiv -M^{\kern 1pt 4}
\!\!\mod 3$, which implies $N\equiv 0\!\!\mod 3$ and $M\equiv 0\!\!\mod 3$. But
$\mu$ and $\nu$ given by the formulas \mythetag{8.11} should be coprime, i\.\,e\.
$\gcd(\mu,\nu)=1$, which contradicts $N\equiv 0\!\!\mod 3$ and $M\equiv 0\!\!\mod 3$. 
The contradiction obtained proves Lemma~\mythelemma{9.1}.\par
     The curve \mythetag{6.1} is associated with the curve \mythetag{6.2}. In the case
of the curve \mythetag{5.3} the associated curve is given by the equation 
$$
\hskip -2em
\tilde y^2=\tilde x^3-15\,\tilde x+22,
\mytag{9.6}
$$
where $\tilde c=-2\,c=2$, $\tilde a=-4\,a-15\,c^2=-15$, and $\tilde B=3\,\tilde c^2
+\tilde a=-3$.\par
\mylemma{9.2} For the curve \mythetag{9.6} the group $\tilde\alpha(\tilde E(\Bbb Q))$ 
has two elements presented by the numbers $1$ and \pagebreak $3$, i\.\,e\. 
$|\tilde\alpha(\tilde E(\Bbb Q))|=2$.
\endproclaim
     The number $\tilde B=-3$ has two positive divisors $1$ and $3$. Therefore, 
according to Lemma~\mythelemma{8.2}, the group $\tilde\alpha(\tilde E(\Bbb Q))$ has at 
most four elements presented by the numbers $1$, $-1$, $3$, $-3$. The formulas
\mythetag{6.19} in this case yield
$$
\xalignat 2
&\hskip -2em
\tilde\alpha(\tilde P_\infty)=1,
&&\tilde\alpha(\tilde P_{\kern 1pt 0})=-3.
\mytag{9.7}
\endxalignat
$$
Since $3=(-1)\cdot(-3)$, it is sufficient to prove that there is no element presented 
by the number $-1$ in $\tilde\alpha(\tilde E(\Bbb Q))$. Assume to the contrary that 
$\tilde\alpha(\tilde P)=-1$ for some rational point $\tilde P$. Due to \mythetag{9.7} 
we have $\tilde P\neq\tilde P_\infty$ and $\tilde P\neq\tilde P_{\kern 1pt 0}$, i\.\,e\. 
$\tilde P=(\tilde x;\,\tilde y)$ is a non-exceptional point. The curve \mythetag{9.6} 
has no rational points with $\tilde y=0$ other than $\tilde P_{\kern 1pt 0}=(2;\,0)$.
Therefore we can write the formulas analogous to \mythetag{8.2}:
$$
\xalignat 2
&\hskip -2em
\tilde x=2+\frac{\tilde m}{\tilde q_1},
&\tilde y=\frac{\tilde n}{\tilde q_2}.
\mytag{9.8}
\endxalignat
$$
Here $\tilde m$, $\tilde n$, $\tilde q_1$, and $\tilde q_2$ are nonzero integer numbers 
such that the fractions \mythetag{9.8} are irreducible. Substituting \mythetag{9.8}
into \mythetag{9.6}, we derive the equality
$$
\hskip -2em
\frac{\tilde n^2}{\tilde q_2^2}=\frac{\tilde m}{\tilde q_1}\cdot
\frac{\tilde m^2+6\,\tilde m\,\tilde q_1-3\,\tilde q_1^2}{\tilde q_1^2} 
\mytag{9.9}
$$
analogous to \mythetag{8.3}. Now, repeating the arguments used in proving
Lemma~\mythelemma{8.1}, we get $\tilde q_1=\tilde q^2$ and $\tilde q_2=\tilde q^3$.
Then we bring \mythetag{9.9} to the equality 
$$
\hskip -2em
\tilde n^2=\tilde m\,(\tilde m^2+6\,\tilde m\,\tilde q^2-3\,\tilde q^4). 
\mytag{9.10}
$$
analogous to \mythetag{8.7}. The positive integer $\tilde\beta$ defined by the formula
$$
\beta=\gcd(|\tilde m|,|\tilde m^2+6\,\tilde m\,\tilde q^2-3\,\tilde q^4|)
=\gcd(|\tilde m|,3) 
$$
is a divisor of the number $\tilde B=-3$. Since above we assumed 
$\tilde\alpha(\tilde P)=-1$, the formula \mythetag{6.18} complemented with 
$\tilde c=2$ yields the formula
$$
\hskip -2em
\tilde\alpha(\tilde P)=\tilde x-2\equiv -1\!\!\!\mod\Bbb Q^{*2}
\mytag{9.11}
$$
analogous to \mythetag{9.2}. From \mythetag{9.11} we derive $\tilde m<0$ and write
$$
\xalignat 2
&\hskip -2em
\tilde m=-\tilde\beta\,\tilde\mu,
&&\tilde m^2+6\,\tilde m\,\tilde q^2-3\,\tilde q^4=-\tilde\beta\,\tilde\nu,
\mytag{9.12}
\endxalignat
$$ 
where $\tilde\mu>0$, $\tilde\nu>0$, and $\gcd(\tilde\mu,\tilde\nu)=1$. The formulas
\mythetag{9.12} are analogous to \mythetag{8.9}. Applying these formulas to 
\mythetag{9.10}, we derive 
$$
\hskip -2em
\tilde n^2=\tilde\beta^{\kern 1pt 2}\kern 0.5pt\tilde\mu\,\tilde\nu. 
\mytag{9.13}
$$
The positive integers $\mu$ and $\nu$ in \mythetag{8.10} are coprime, while their 
product is a square. Therefore both of these integer numbers are squares:
$$
\xalignat 2
&\hskip -2em
\tilde\mu=\tilde M^{\kern 1pt 2},
&&\tilde\nu=\tilde N^{\kern 1pt 2}.
\mytag{9.14}
\endxalignat
$$ 
The formulas \mythetag{9.13} and \mythetag{9.14} are analogous to \mythetag{8.10} and 
\mythetag{8.11} \pagebreak respectively. Applying \mythetag{9.14} to \mythetag{9.12} 
and then applying \mythetag{9.12} to \mythetag{9.8}, we derive
$$
\hskip -2em
\tilde x-2=\frac{\tilde m}{\tilde q^2}=-\frac{\tilde\beta\,\tilde\mu}{\tilde q^2}
=-\frac{\tilde\beta\,\tilde M^{\kern 1pt 2}}{\tilde q^2}\equiv-\tilde \beta\!\!\!
\mod\Bbb Q^{*2},
\mytag{9.15}
$$
Comparing \mythetag{9.15} with \mythetag{9.11}, we conclude that $\tilde\beta$ is a
square, i\.\,e\. $\tilde\beta=\tilde K^{\kern 1pt 2}$ for some nonzero integer
number $\tilde K$. On the other hand, in the present case $\tilde\beta$ is a divisor 
of $\tilde B=-3$. Therefore $\tilde\beta=\tilde K^{\kern 1pt 2}=1$. As a result from
\mythetag{9.12} and \mythetag{9.13}, we derive 
$$
\xalignat 2
&\hskip -2em
\tilde m=-\tilde M^{\kern 1pt 2},
&&\tilde n^2=\tilde M^{\kern 1pt 2}\,\tilde N^{\kern 1pt 2}.
\mytag{9.16}
\endxalignat
$$
Substituting \mythetag{9.16} into \mythetag{9.10}, we produce the equality 
analogous to \mythetag{9.5}:
$$
\hskip -2em
\tilde N^{\kern 1pt 2}=-\tilde M^{\kern 1pt 4}+6\,\tilde M^{\kern 1pt 2}
\,\tilde q^2+3\,\tilde q^4. 
\mytag{9.17}
$$
The equality \mythetag{9.17} means that $\tilde N^{\kern 1pt 2}\equiv 
-\tilde M^{\kern 1pt 4}\!\!\mod 3$, which implies $\tilde N\equiv 0\!\!\mod 3$ 
and $\tilde M\equiv 0\!\!\mod 3$. But $\tilde\mu$ and $\tilde\nu$ given by the 
formulas \mythetag{9.14} should be coprime, i\.\,e\. $\gcd(\tilde \mu,\tilde \nu)=1$, 
which contradicts $\tilde N\equiv 0\!\!\mod 3$ and $\tilde M\equiv 0\!\!\mod 3$. 
The contradiction obtained proves Lemma~\mythelemma{9.2}.\par
      Let's apply Lemma~\mythelemma{7.1} to the curve \mythetag{5.3}. Substituting
$y=0$ into \mythetag{5.3}, we derive $x^3+1=0$. The left hand side of this 
equation factors as follows:
$$
\hskip -2em
(x+1)\,(x^2-x+1)=0.
\mytag{9.18}
$$
It is clear that $x=-1$ is the only rational solution of the equation 
\mythetag{9.18}. If we recall that $c=-1$ for the curve \mythetag{5.3} and 
$P_{\kern 1pt 0}=(c;\,0)$, then we derive 
$$
\hskip -2em
[\Ker\tildepsi:(\Ker\tildepsi\cap\psi(E(\Bbb Q))]=2. 
\mytag{9.19}
$$
Let's substitute \mythetag{9.19} into the formula \mythetag{7.21}. Then, applying 
Lemma~\mythelemma{9.1} and Lemma~\mythelemma{9.2} to the formula \mythetag{7.21},
we obtain
$$
\hskip -2em
2^{\kern 1pt r}\!\cdot|E_2(\Bbb Q)|=2.
\mytag{9.20}
$$
The subgroup $E_2(\Bbb Q)\subset E_{\sssize\text{tors}}(\Bbb Q)$ in \mythetag{9.20}
and the torsion group $E_{\sssize\text{tors}}(\Bbb Q)$ itself can be calculated with 
the use of the following Lutz-Nagell theorem.\par
\mytheoremwithtitle{9.1}{ (Lutz {\rm and} Nagell)} Let $P=(x;\, y)$ be  a rational 
point of finite order on an elliptic curve given by the equation \mythetag{2.3} 
with integer coefficients\/ $a$ and $b$. Then its coordinates $x$ and $y$ 
both are integers and either $y=0$ or $y$ divides $D$, where 
$D=-4\,a^3-27\,b^{\kern 1pt 2}$ is the discriminant of the cubic polynomial in the 
right hand side of the equation \mythetag{2.3}.
\endproclaim
     The discriminant of the cubic polynomial $x^3+1$ in the case of the curve 
\mythetag{5.3} is calculated explicitly: $D=-27$. Then Theorem~\mythetheorem{9.1}
yields 
$$
E_{\sssize\text{tors}}(\Bbb Q)=\{(\infty;\,\infty),\,(2;\,3),\,(0;\,1),\,(-1;\,0),\,
(0;\,-1),\,(2;\,-3)\}\cong\Bbb Z_6.
$$
The group $\Bbb Z_6$ has one element $P_{\kern 1pt 0}=(-1;\,0)$ of the order $2$.
Therefore  
$$
\pagebreak
E_2(\Bbb Q)=\{(\infty;\,\infty),\,(-1;\,0)\}\cong\Bbb Z_2.
$$
and $|E_2(\Bbb Q)|=2$. Applying this result to \mythetag{9.20}, we derive $2^{\kern 
1pt r}=1$, which means $r=0$. The proof of Theorem~\mythetheorem{5.6} is 
completed. As we noted above in section 5, this theorem is equivalent to Euler's 
theorem~\mythetheorem{5.5}.\par
     Theorem~\mythetheorem{5.6} is very important in the context of cuboid curves 
\mythetag{2.2}. As we have seen in section 5, it is used for proving 
Theorem~\mythetheorem{5.7}, which solves the rank problem for the subset of curves 
\mythetag{2.2} where $4\,R^{\kern 0.5pt 2}\,N$ is an exact cube. The subset of curves 
where $4\,R^{\kern 0.5pt 2}\,N$ is not an exact cube is much more broad and more 
complicated. Our further efforts below are directed toward solving the rank problem 
for this subset of curves \mythetag{2.2}, though we do not solve this problem 
completely.\par
\head
10. Three descent. Isogenies and descent mappings.
\endhead
      From this section on, below we assume that the integer number 
$4\,R^{\kern 0.5pt 2}\,N$ is not an exact cube. Then, applying the transformation 
\mythetag{5.2} to the curve equation \mythetag{2.2}, we can bring this curve 
equation to the form 
$$
\hskip -2em
y^2=x^3+e^2,
\mytag{10.1}
$$
where $e\neq 1$ is some positive cube free integer number. Like in the case of the 
curves \mythetag{6.1} and \mythetag{6.2}, along with each curve $E$ of the form 
\mythetag{10.1}, we consider its associated curve $\tilde E$. This curve is defined
by means of the equation
$$
\hskip -2em
\tilde y^2=\tilde x^3-27\,e^2.
\mytag{10.2}
$$
Applying Theorem~\mythetheorem{5.3} to the curve \mythetag{10.1}, we derive the 
following lemma.
\mylemma{10.1} If $e\neq 1$ is a positive third power free integer, then the 
torsion subgroup $E_{\sssize\text{tors}}(\Bbb Q)$ of the curve \mythetag{10.1} 
is isomorphic to $\Bbb Z_3$:
$$
\hskip -2em
E_{\sssize\text{tors}}(\Bbb Q)=\{(\infty;\,\infty),\,(0;\,e),\,(0;\,-e)\}\cong\Bbb Z_3.
\mytag{10.3}
$$
\endproclaim
     The points $P^{\kern 0.5pt \sssize +}_0=(0;\,e)$ and 
$P^{\kern 0.5pt \sssize -}_0=(0;\,-e)$ in \mythetag{10.3} correspond to the exceptional
points \mythetag{4.1} of the initial elliptic curve \mythetag{2.2}.\par
     The associated curve $\tilde E$ in \mythetag{10.2} is different. Applying 
Theorem~\mythetheorem{5.3}, we derive the following lemma for the associated curve 
\mythetag{10.2}.
\mylemma{10.2} If $e=4$, then the torsion subgroup $\tilde E_{\sssize\text{tors}}(\Bbb Q)$ 
of the curve \mythetag{10.2} is isomorphic to $\Bbb Z_3$. It is given explicitly by the 
formula
$$
\tilde E_{\sssize\text{tors}}(\Bbb Q)=\{(\infty;\,\infty),\,(12,\,36),\,(12,\,-36)\},
$$
where $P^{\kern 0.5pt \sssize +}_0=(12;\,36)$ and $P^{\kern 0.5pt \sssize -}_0=(12;\,-36)$
are two third order rational points of the curve $\tilde E$. If\/ $e\neq 1$ and $e\neq 4$
is a positive third power free integer, then the torsion subgroup 
$\tilde E_{\sssize\text{tors}}(\Bbb Q)$ is trivial, i\.\,e\. 
$\tilde E_{\sssize\text{tors}}(\Bbb Q)=\{(\infty;\,\infty)\}$.
\endproclaim
      Like in Section 6, we define a mapping $\psi\!:\,E\longrightarrow\tilde E$ by setting
$$
\xalignat 2
&\hskip -2em
\tilde x=\frac{x^3+4\,e^2}{x^2},
&&\tilde y=\frac{y\,(x^3-8\,e^2)}{x^3},
\mytag{10.4}
\endxalignat
$$
where $x\neq 0$. For the exceptional points $P^{\kern 0.5pt \sssize +}_0$,
$P^{\kern 0.5pt \sssize -}_0$, and $P_\infty$ we set by definition
$$
\pagebreak
\xalignat 3
&\hskip -2em
\psi(P^{\kern 0.5pt \sssize +}_0)=\tilde P_\infty,
&&\psi(P^{\kern 0.5pt \sssize -}_0)=\tilde P_\infty,
&&\psi(P_\infty)=\tilde P_\infty.\quad
\mytag{10.5}
\endxalignat
$$
The formulas \mythetag{10.4}  are analogs of the formulas \mythetag{6.3} and \mythetag{6.4},
while the formulas \mythetag{10.5} are analogs of the formulas \mythetag{6.5}.\par
     The curve \mythetag{10.2} is similar to the curve \mythetag{10.1}. Therefore there is a
backward mapping $\tildepsi\!:\,\tilde E\longrightarrow E$. This mapping is given by the 
formulas  
$$
\xalignat 2
&\hskip -2em
x=\frac{\tilde x^3-108\,e^2}{9\,\tilde x^2},
&&y=\frac{\tilde y\,(\tilde x^3+216\,e^2)}{27\,\tilde x^3},
\mytag{10.6}
\endxalignat
$$
where $\tilde x\neq 0$. The formulas \mythetag{10.6} are analogs of the formulas 
\mythetag{6.6} and \mythetag{6.7}. The curve \mythetag{10.2} has no exceptional rational 
points with $\tilde x=0$. Therefore the analog of the formulas \mythetag{6.8} in this
case is written as follows:
$$
\tildepsi(\tilde P_\infty)=P_\infty.\quad
\mytag{10.7}
$$
The transformations \mythetag{10.4} and \mythetag{10.6} were discovered by Fueter
in \mycite{12} (see also \mycite{13} and Chapter~26 of \mycite{14}).\par
\mylemma{10.3} The mapping $\psi\!:\,E\longrightarrow\tilde E$ defined by the formulas
\mythetag{10.4} and \mythetag{10.5} induces a homomorphism of Abelian groups
$\psi\!:\,E(\Bbb Q)\longrightarrow\tilde E(\Bbb Q)$. 
\endproclaim
\mylemma{10.4} The mapping $\tildepsi\!:\,\tilde E\longrightarrow E$ defined by the 
formulas \mythetag{10.6} and \mythetag{10.7} induces the homomorphism of Abelian groups 
$\tildepsi\!:\,\tilde E(\Bbb Q)\longrightarrow E(\Bbb Q)$. 
\endproclaim
      Lemmas~\mythelemma{10.3} and \mythelemma{10.4} are analogs of Lemmas~\mythelemma{6.2} 
and \mythelemma{6.3} respectively. Their proofs are pure
calculations using the formulas \mythetag{3.3}, \mythetag{3.6}, \mythetag{3.7}, 
\mythetag{3.8}, \mythetag{3.9}.\par
     Now let's consider the composite mapping $\tildepsi\compos\psi\!:\,E\longrightarrow E$.
Due to Lemmas~\mythelemma{10.3} and \mythelemma{10.4} it induces an endomorphism of
the Abelian group $E(\Bbb Q)$. 
\mylemma{10.5} The endomorphism $\tildepsi\compos\psi\!:\,E(\Bbb Q)\longrightarrow 
E(\Bbb Q)$ coincides with the tripling endomorphism, i\.\,e\. $\tildepsi\compos\psi(P)
=3\,P$ for any $P\in E(\Bbb Q)$.
\endproclaim
     The proof of Lemma~\mythelemma{10.5} is also pure calculations.\par 
     Recall that $\Bbb Q^*$ is the set of all nonzero rational numbers. This set possesses 
the structure of a multiplicative Abelian group. Through $\Bbb Q^{*3}$ we denote the set
of all nonzero rational numbers which are cubes. Then $\Bbb Q^{*3}$ is a subgroup of 
$\Bbb Q^*$ and we have the factor group $\Bbb Q^*\!/\Bbb Q^{*3}$. Now we define two mappings 
$\alpha_{\sssize +}\!:\,E\longrightarrow\Bbb Q^*\!/\Bbb Q^{*3}$ and
$\alpha_{\sssize -}\!:\,E\longrightarrow\Bbb Q^*\!/\Bbb Q^{*3}$. The mapping 
$\alpha_{\sssize -}$ is defined by means of the formula  
$$
\hskip -2em
\alpha_{\sssize -}(P)=y-e\text{\ \ for \ }y\neq e\text{\ \ and \ }y\neq\infty.
\mytag{10.8}
$$
The exceptional points $P^{\kern 0.5pt \sssize +}_0$ and $P_\infty$ are 
treated separately. For them we set
$$
\xalignat 2
&\hskip -2em
\alpha_{\sssize -}(P_\infty)=1,
&&\alpha_{\sssize -}(P^{\kern 0.5pt \sssize +}_0)=\frac{1}{2\,e}.
\mytag{10.9}
\endxalignat
$$
The mapping $\alpha_{\sssize +}$ is defined similarly by means of the formula  
$$
\hskip -2em
\alpha_{\sssize +}(P)=y+e\text{\ \ for \ }y\neq -e\text{\ \ and \ }y\neq\infty.
\mytag{10.10}
$$
The exceptional points in this case are $P^{\kern 0.5pt \sssize -}_0$ and 
$P_\infty$. For them we set
$$
\xalignat 2
&\hskip -2em
\alpha_{\sssize +}(P_\infty)=1,
&&\alpha_{\sssize +}(P^{\kern 0.5pt \sssize -}_0)=-\frac{1}{2\,e}.
\mytag{10.11}
\endxalignat
$$
The formulas \mythetag{10.8} and \mythetag{10.9} are taken from \mycite{15}. The 
formulas \mythetag{10.10} and \mythetag{10.11} are written by analogy.\par
      In order to relate \mythetag{10.9} with \mythetag{10.8} and in order to
relate \mythetag{10.11} with \mythetag{10.10} we write the equation \mythetag{10.1} 
as $(y+e)(y-e)=x^3$. Then \mythetag{10.8} and \mythetag{10.10} yield 
$$
\xalignat 2
&\hskip -2em
\alpha_{\sssize -}(P)=\frac{x^3}{y+e}, 
&&\alpha_{\sssize +}(P)=\frac{x^3}{y-e}, 
\mytag{10.12}
\endxalignat
$$
Cubic factors are neglected modulo $\Bbb Q^{*3}$ in $\Bbb Q^*\!/\Bbb Q^{*3}$. 
Therefore, the above formulas \mythetag{10.12} are equivalent to the following 
four relationships:
$$
\xalignat 2
&\hskip -2em
\alpha_{\sssize -}(P)=\frac{1}{y+e}, 
&&\alpha_{\sssize +}(P)=\frac{1}{y-e}, 
\mytag{10.13}\\
\vspace{1ex}
&\hskip -2em
\alpha_{\sssize -}(P)=\frac{y^3}{x^3(\kern 1pt y+e)}, 
&&\alpha_{\sssize +}(P)=\frac{y^3}{x^3(\kern 1pt y-e)}.
\mytag{10.14}\\
\endxalignat
$$
Passing to the limit as $y\to e$ and as $y\to -e$ in \mythetag{10.13} and passing
to the limit as $x\to\infty$ along with $y\to\infty$ in \mythetag{10.14}, we 
derive the formulas \mythetag{10.9} and \mythetag{10.11}. 
\mylemma{10.6} The mapping $\alpha_{\sssize -}\!:\,E\longrightarrow
\Bbb Q^*\!/\Bbb Q^{*3}$ defined by the formulas \mythetag{10.8} and \mythetag{10.9} 
induces the homomorphism of Abelian groups $\alpha_{\sssize -}\!:\,E(\Bbb Q)
\longrightarrow\Bbb Q^*\!/\Bbb Q^{*3}$. 
\endproclaim
\mylemma{10.7} The mapping $\alpha_{\sssize +}\!:\,E\longrightarrow
\Bbb Q^*\!/\Bbb Q^{*3}$ defined by the formulas \mythetag{10.10} and \mythetag{10.11} 
induces the homomorphism of Abelian groups $\alpha_{\sssize +}\!:\,E(\Bbb Q)
\longrightarrow\Bbb Q^*\!/\Bbb Q^{*3}$. 
\endproclaim
     In order to prove Lemma~\mythelemma{10.6} assume that $P_{\kern 1pt 1}
=(x_1;\,y_1)$ and $P_{\kern 1pt 2}=(x_2;\,y_2)$ are two non-exceptional rational 
points of the curve \mythetag{10.1} and assume that the point 
$P_{\kern 1pt 3}=(x_3;\,y_3)$ is their sum. Then \mythetag{10.8} yields 
$\alpha(P_{\kern 1pt 3})=y_3-e$. Applying the formula \mythetag{3.7}, where $a=0$ 
and $b=e^2$ in the case of the curve \mythetag{10.1}, we derive
$$
\gathered
\alpha_{\sssize -}(P_{\kern 1pt 3})
=(\kern 1pt y_2\,x_1^3-y_1\,x_2^3+3\,y_2\,x_2\,x_1^2-3\,y_1\,x_1\,x_2^2
+4\,y_2\,e^2-4\,y_1\,e^2+\\
+\,e\,x_2^3-e\,x_1^3+3\,e\,x_2\,x_1^2-3\,e\,x_1\,x_2^2)\,(x_1-x_2)^{-3}.
\endgathered\quad
\mytag{10.15}
$$ 
Now let's calculate the product $\alpha_{\sssize -}(P_{\kern 1pt 1})
\,\alpha_{\sssize -}(P_{\kern 1pt 2})$. Applying \mythetag{10.8}, we get 
$$
\hskip -2em
\alpha_{\sssize -}(P_{\kern 1pt 1})\ \alpha_{\sssize -}(P_{\kern 1pt 2})
=(\kern 1pt y_1-e)\,(\kern 1pt y_2-e)=\frac{(\kern 1pt y_1^2-e^2)
\,(\kern 1pt y_2^2-e^2)}{(\kern 1pt y_1+e)\,(\kern 1pt y_2+e)}.
\mytag{10.16}
$$
Then we apply the curve equation to the numerator of the fraction in 
\mythetag{10.16}:
$$
\hskip -2em
\alpha_{\sssize -}(P_{\kern 1pt 1})\ \alpha_{\sssize -}(P_{\kern 1pt 2})
=\frac{x_1^3\,x_2^3}
{(\kern 1pt y_1+e)\,(\kern 1pt y_2+e)}.
\mytag{10.17}
$$
Cubic factors are neglected modulo $\Bbb Q^{*3}$ in $\Bbb Q^*\!/\Bbb Q^{*3}$.
Keeping in mind this fact, from \mythetag{10.15} and \mythetag{10.17} 
we derive the following relationship:
$$
\gathered
\alpha_{\sssize -}(P_{\kern 1pt 3})\ \alpha_{\sssize -}(P_{\kern 1pt 1})^{-1}
\alpha_{\sssize -}(P_{\kern 1pt 2})^{-1}=(\kern 1pt y_2\,x_1^3-y_1\,x_2^3+3\,y_2\,x_2\,x_1^2-3\,y_1\,x_1\,x_2^2\,+\\
+\,4\,y_2\,e^2-4\,y_1\,e^2+e\,x_2^3-e\,x_1^3+3\,e\,x_2\,x_1^2
-3\,e\,x_1\,x_2^2)\,{(\kern 1pt y_1+e)\,(\kern 1pt y_2+e)}.
\endgathered\quad
\mytag{10.18}
$$
Taking into account the curve equation \mythetag{10.1}, \pagebreak we can 
expand and then refactor the right hand side of the equality \mythetag{10.18}. 
As a result we get 
$$
\gathered
\alpha_{\sssize -}(P_{\kern 1pt 3})\ \alpha_{\sssize -}(P_{\kern 1pt 1})^{-1}
\alpha_{\sssize -}(P_{\kern 1pt 2})^{-1}=(e\,x_2-e\,x_1+y_1\,x_2-x_1\,y_2)^3.
\endgathered\quad
\mytag{10.19}
$$
Cubic factors are neglected modulo $\Bbb Q^{*3}$ in $\Bbb Q^*\!/\Bbb Q^{*3}$.
Therefore \mythetag{10.19} yields 
$$
\hskip -2em
\alpha_{\sssize -}(P_{\kern 1pt 3})=\alpha_{\sssize -}(P_{\kern 1pt 1})
\ \alpha_{\sssize -}(P_{\kern 1pt 2}).
\mytag{10.20}
$$
The equality \mythetag{10.20} proves Lemma~\mythelemma{10.6}. 
Lemma~\mythelemma{10.7} is proved similarly.\par
     The homomorphisms $\alpha_{\sssize +}\!:\,E\longrightarrow
\Bbb Q^*\!/\Bbb Q^{*3}$ and $\alpha_{\sssize -}\!:\,E\longrightarrow
\Bbb Q^*\!/\Bbb Q^{*3}$ are not independent. Indeed, due to \mythetag{10.8}, 
\mythetag{10.10}, and the curve equation \mythetag{10.1} we have
$\alpha_{\sssize +}(P)\,\alpha_{\sssize -}(P)=x^3$. But cubic factors are 
neglected modulo $\Bbb Q^{*3}$ in $\Bbb Q^*\!/\Bbb Q^{*3}$. Then
$$
\hskip -2em
\alpha_{\sssize +}(P)\,\alpha_{\sssize -}(P)=1.
\mytag{10.21}
$$
The formula \mythetag{10.21} means that Lemma~\mythelemma{10.7} is immediate 
from Lemma~\mythelemma{10.6}. Moreover, from the formula \mythetag{10.21} we 
derive the equality
$$
\hskip -2em
\Ker\alpha_{\sssize -}=\Ker\alpha_{\sssize +}.
\mytag{10.22}
$$\par
\mylemma{10.8} The kernel of the homomorphism $\alpha_{\sssize -}\!:\,E(\Bbb Q)
\longrightarrow\Bbb Q^*\!/\Bbb Q^{*3}$ coincides with the image of the
homomorphism $\tildepsi\!:\,\tilde E(\Bbb Q)\longrightarrow E(\Bbb Q)$.
\endproclaim
\mylemma{10.9} The kernel of the homomorphism $\alpha_{\sssize +}\!:\,E(\Bbb Q)
\longrightarrow\Bbb Q^*\!/\Bbb Q^{*3}$ coincides with the image of the
homomorphism $\tildepsi\!:\,\tilde E(\Bbb Q)\longrightarrow E(\Bbb Q)$.
\endproclaim
     Due to the formula \mythetag{10.21} Lemma~\mythelemma{10.9} is immediate 
from Lemma~\mythelemma{10.8}. Therefore it is sufficient to prove  
Lemma~\mythelemma{10.8}.\par
     Let $P=(x;\,y)$ be a non-exceptional rational point of the curve \mythetag{10.1}
such that $P\in\Img\tildepsi$. Then its coordinates $x$ and $y$ are given by the
formulas \mythetag{10.6}. Using the second formula \mythetag{10.6} and applying
\mythetag{10.8}, we derive the following relationship:
$$
\hskip -2em
\alpha_{\sssize -}(P)=y-e=\frac{216\,\tilde y\,e^2+\tilde y\,\tilde x^3
-27\,e\,\tilde x^3}{27\,\tilde x^3}.
\mytag{10.23}
$$
Due to the curve equation \mythetag{10.2} the numerator of the fraction in the right 
hand side of \mythetag{10.23} factors as $216\,\tilde y\,e^2+\tilde y\,\tilde x^3
-27\,e\,\tilde x^3=(\tilde y-9\,e)^3$. Then
$$
\hskip -2em
\alpha_{\sssize -}(P)=\frac{(\tilde y-9\,e)^3}{3^3\,\tilde x^3}.
\mytag{10.24}
$$
Cubic factors are neglected modulo $\Bbb Q^{*3}$ in $\Bbb Q^*\!/\Bbb Q^{*3}$. 
Therefore \mythetag{10.24} means that $\alpha_{\sssize -}(P)=1$ and 
$P\in\Ker\alpha_{\sssize -}$.\par
     For the exceptional point $P=P^{\kern 0.5pt \sssize -}_0=(0;\,-e)$ the value 
of $\alpha_{\sssize -}(P)$ is given by the formula \mythetag{10.8}, which serves
non-exceptional points as well. Therefore the above reasons remain valid for the
exceptional point $P=P^{\kern 0.5pt \sssize -}_0$, i\.\,e\. if $P^{\kern 0.5pt 
\sssize -}_0\in\Img\tildepsi$, then $P^{\kern 0.5pt \sssize -}_0\in\Ker
\alpha_{\sssize -}$. The exceptional point $P=P^{\kern 0.5pt \sssize +}_0=(0;\,e)$
is more special. Assume that $P^{\kern 0.5pt \sssize +}_0\in\Img\tildepsi$. Then
$P^{\kern 0.5pt \sssize +}_0=\tildepsi(\tilde P)$ for some non-exceptional 	point
$\tilde P=(\tilde x;\,\tilde y)$ of the curve \mythetag{10.2}. Applying the formulas
\mythetag{10.6}, we derive the following two equations: 
$$
\xalignat 2
&\hskip -2em
\frac{\tilde x^3-108\,e^2}{9\,\tilde x^2}=0,
&&\frac{\tilde y\,(\tilde x^3+216\,e^2)}{27\,\tilde x^3}=e.
\mytag{10.25}
\endxalignat
$$
Since $\tilde x\neq 0$, the equations \mythetag{10.25} are equivalent to
$$
\xalignat 2
&\hskip -2em
\tilde x^3-108\,e^2=0,
&&216\,\tilde y\,e^2+\tilde y\,\tilde x^3-27\,e\,\tilde x^3=0.
\quad
\mytag{10.26}
\endxalignat
$$
The left hand side of the second equation \mythetag{10.26} coincides with the numerator
of the fraction in \mythetag{10.23}. Therefore, using the curve equation 
\mythetag{10.2}, we can simplify it to the equation $(\tilde y-9\,e)^3=0$, which
yields $\tilde y=9\,e$. Substituting $\tilde y=9\,e$ into the curve equation
\mythetag{10.2}, we derive the equation $\tilde x^3=108\,e^2$, which coincides with
the first equation \mythetag{10.25}.\par
     Let's recall that $e$ is a third power free positive integer 
(see \mythetag{10.1}) and use the prime factor expansion $108=2^2\kern -1pt
\cdot 3^3$. Then the equality $\tilde x^3=108\,e^2$ can hold in the only case where 
$e=4$. For the coordinates of the point $\tilde P$ we derive $\tilde x=12$ and 
$\tilde y=36$. The image of the point $\tilde P=(12;\,36)$ under the mapping 
$\tildepsi$ is the exceptional point $P^{\kern 0.5pt \sssize +}_0=(0;\,4)$. 
The value of $\alpha_{\sssize -}(P^{\kern 0.5pt \sssize +}_0)$ for such a point is 
given by the formula \mythetag{10.9}. Since $e=4$, this formula yields 
$\alpha_{\sssize -}(P^{\kern 0.5pt \sssize +}_0)=1/8=(1/2)^3$. Cubic factors are 
neglected modulo $\Bbb Q^{*3}$ in $\Bbb Q^*\!/\Bbb Q^{*3}$. Therefore 
$\alpha_{\sssize -}(P^{\kern 0.5pt \sssize +}_0)=1/8$ is equivalent to 
$\alpha_{\sssize -}(P^{\kern 0.5pt \sssize +}_0)=1$, which means 
$P^{\kern 0.5pt \sssize +}_0\in\Ker\alpha_{\sssize -}$.\par 
     Thus we have proved the inclusion  $\Img\tildepsi
\subseteq\Ker\alpha_{\sssize -}$.\par
     Now assume that $P=(x;\,y)$ is some non-exceptional rational point of 
the curve \mythetag{10.1} such that $P\in\Ker\alpha_{\sssize -}$. Then 
$\alpha_{\sssize -}(P)=1$, which means 
$$
\hskip -2em
y-e=z^3\text{, \ where \ }z\neq 0\text{\ \ and \ }z\in\Bbb Q.
\mytag{10.27}
$$
The formula \mythetag{10.27} yields $y=z^3+e$. Substituting $y=z^3+e$ into the curve
equation \mythetag{10.1}, we derive the following equation for $x$ and $z$: 
$$
\hskip -2em
z^6+2\,z^3\,e-x^3=0.
\mytag{10.28}
$$
Using $x$ and $z$, we define the point $\tilde P=(\tilde x;\,\tilde y)$ of the 
associated curve \mythetag{10.2} whose coordinates $\tilde x$ and $\tilde y$ are given 
by the formulas 
$$
\xalignat 2
&\hskip -2em
\tilde x=\frac{6\,e\,z}{x-z^2},
&&\tilde y=\frac{9\,e\,(x+z^2)}{x-z^2}.
\mytag{10.29}
\endxalignat
$$
The denominator of the fractions \mythetag{10.29} is nonzero. Indeed, otherwise we would
have $x=z^2$. Substituting $x=z^2$ into \mythetag{10.28}, we would derive $2\,z^3\,e=0$.
Since $e\neq 0$, this would mean $z=0$, which contradicts \mythetag{10.27}.\par
     In order to verify that $\tilde x$ and $\tilde y$ given by the formulas 
\mythetag{10.29} do actually form a point of the curve \mythetag{10.2}, we substitute them
into the equation \mythetag{10.2}. This yields
$$
\hskip -2em
\frac{108\,e^2\,(z^6+2\,z^3\,e-x^3)}{(x-z^2)^3}=0.
\mytag{10.30}
$$
It is easy to see that the equation \mythetag{10.30} is fulfilled due to \mythetag{10.28}. 
\par     
     The homomorphism $\tildepsi\!:\,\tilde E(\Bbb Q)\longrightarrow E(\Bbb Q)$ is defined 
by means of the formulas \mythetag{10.6}. Substituting \mythetag{10.29} into \mythetag{10.6}, 
we derive 
$$
\allowdisplaybreaks
\gather
\hskip -2em
x=\frac{2\,z^3\,e-x^3+3\,x^2\,z^2-3\,x\,z^4+z^6}{3\,z^2\,(x-z^2)},
\mytag{10.31}\\
\vspace{1ex}
\hskip -2em
y=\frac{(z^3\,e+x^3-3\,x^2\,z^2+3\,x\,z^4-z^6)\,(x+z^2)}{3\,z^3\,(x-z^2)}.
\mytag{10.32}
\endgather
$$
The polynomials in the numerators of the fractions \mythetag{10.31} and \mythetag{10.32} 
are simplified with the use of the equation \mythetag{10.28}. As a result the equality 
\mythetag{10.31} turns to the identity $x=x$. The equality \mythetag{10.32} turns to
the equality $y=e+z^3$, which is equivalent to \mythetag{10.27}. Thus we have proved that
$P=\tildepsi(\tilde P)$, where $\tilde P$ is the point of the curve \mythetag{10.2} with
the coordinates \mythetag{10.29}, i\.\,e\. $P\in\Img\tildepsi$.\par
     The exceptional point $P=P^{\kern 0.5pt \sssize -}_0=(0;\,-e)$ is similar to 
non-exceptional points. The value of $\alpha_{\sssize -}(P)$ for it is given by the
formula \mythetag{10.8}, which serves non-exceptio\-nal points as well. Therefore the above 
reasons remain valid for the exceptional point $P=P^{\kern 0.5pt \sssize -}_0$, i\.\,e\. 
if $P^{\kern 0.5pt \sssize -}_0\in\Ker\alpha_{\sssize -}$, then $P^{\kern 0.5pt \sssize -}_0
\in\Img\tildepsi$. The other exceptional point $P=P^{\kern 0.5pt \sssize +}_0=(0;\,e)$
is more special. Assume that $P^{\kern 0.5pt \sssize +}_0\in\Ker\alpha_{\sssize -}$. 
Then $\alpha_{\sssize -}(P^{\kern 0.5pt \sssize +}_0)=1$. Applying the second formula
\mythetag{10.9}, we derive the following equation: 
$$
\hskip -2em
\frac{1}{2\,e}=z^3\text{, \ where \ }z\neq 0\text{\ \ and \ }z\in\Bbb Q.
\mytag{10.33}
$$
Let's recall that $e$ is a third power free positive integer (see \mythetag{10.1}).
Then the equality \mythetag{10.33} can hold in the only case where $e=4$. In this case
$P^{\kern 0.5pt \sssize +}_0=(0;\,4)$. As we have seen above, if $e=4$ the point 
$P^{\kern 0.5pt \sssize +}_0=(0;\,4)$ is the image of the point $\tilde P=(12;\,36)$ under 
the mapping $\tildepsi$, i\.\,e\. $P^{\kern 0.5pt \sssize +}_0\in\Img\tildepsi$.\par
     Thus we have proved the inclusion $\Ker\alpha_{\sssize -}\subseteq\Img\tildepsi$.
Combining it with the previously proved inclusion $\Img\tildepsi\subseteq\Ker
\alpha_{\sssize -}$, we derive the equality $\Ker\alpha_{\sssize -}=\Img\tildepsi$, which
proves Lemma~\mythelemma{10.8}.\par
     The curve \mythetag{10.2} is different from the curve \mythetag{10.1}. In order to serve 
this curve we introduce the number field $\Bbb K=\Bbb Q(\sqrt{-3\,})$. Elements of
$\Bbb K$ are presented as 
$$
\hskip -2em
z=u+v\,\sqrt{-3\,}\text{, \ where \ }u\in\Bbb Q \text{\ \ and \ }v\in\Bbb Q.
\mytag{10.34}
$$
Nonzero numbers of the form \mythetag{10.34} constitute a multiplicative Abelian group.
We denote it $\Bbb K^*$ and consider its factor group $\Bbb K^*\!/\Bbb K^{*3}$. Then 
we define two mappings $\tilde\alpha_{\sssize +}\!:\,\tilde E\longrightarrow\Bbb K^*\!
/\Bbb K^{*3}$ and $\alpha_{\sssize -}\!:\,\tilde E\longrightarrow\Bbb K^*\!/\Bbb K^{*3}$. 
The mapping $\tilde\alpha_{\sssize -}$ is defined by means of the formula similar
to the formula \mythetag{10.8}:
$$
\hskip -2em
\tilde\alpha_{\sssize -}(\tilde P)=\tilde y-3\,\sqrt{-3\,}\,e\text{\ \ for \ }
y\neq\infty.
\mytag{10.35}
$$
The exceptional point at infinity $\tilde P_\infty$ is treated separately. For this
point we set
$$
\hskip -2em
\tilde\alpha_{\sssize -}(\tilde P_\infty)=1.
\mytag{10.36}
$$
The mapping $\tilde\alpha_{\sssize +}$ is defined similarly by means of the formula  
$$
\hskip -2em
\tilde\alpha_{\sssize +}(\tilde P)=\tilde y+3\,\sqrt{-3\,}\,e\text{\ \ for \ }
y\neq\infty.
\mytag{10.37}
$$
For the exceptional point at infinity $\tilde P_\infty$ in this case we again set
$$
\hskip -2em
\tilde\alpha_{\sssize +}(\tilde P_\infty)=1.
\mytag{10.38}
$$
The formulas \mythetag{10.35}, \mythetag{10.36}, \mythetag{10.37}, \mythetag{10.38}
\pagebreak are written by analogy to the formulas \mythetag{10.8}, \mythetag{10.9}, 
\mythetag{10.10}, \mythetag{10.11}.\par
\mylemma{10.10} The mapping $\tilde\alpha_{\sssize -}\!:\,\tilde E\longrightarrow
\Bbb K^*\!/\Bbb K^{*3}$ defined by the formulas \mythetag{10.35} and \mythetag{10.36} 
induces the homomorphism of Abelian groups $\tilde\alpha_{\sssize -}\!:\,\tilde E(\Bbb Q)
\longrightarrow\Bbb K^*\!/\Bbb K^{*3}$. 
\endproclaim
\mylemma{10.11} The mapping $\tilde\alpha_{\sssize +}\!:\,\tilde E\longrightarrow
\Bbb K^*\!/\Bbb K^{*3}$ defined by the formulas \mythetag{10.37} and \mythetag{10.38} 
induces the homomorphism of Abelian groups $\tilde\alpha_{\sssize +}\!:\,\tilde E(\Bbb Q)
\longrightarrow\Bbb K^*\!/\Bbb K^{*3}$. 
\endproclaim
     Lemmas~\mythelemma{10.10} and \mythelemma{10.11} are similar to Lemmas~\mythelemma{10.6} 
and \mythelemma{10.7}. They are proved by analogy to Lemmas~\mythelemma{10.6} 
and \mythelemma{10.7}. 
\mylemma{10.12} The kernel of the homomorphism $\tilde\alpha_{\sssize -}\!:\,\tilde E(\Bbb Q)
\longrightarrow\Bbb K^*\!/\Bbb K^{*3}$ coincides with the image of the
homomorphism $\psi\!:\,E(\Bbb Q)\longrightarrow\tilde E(\Bbb Q)$.
\endproclaim
\mylemma{10.13} The kernel of the homomorphism $\tilde\alpha_{\sssize +}\!:\,\tilde E(\Bbb Q)
\longrightarrow\Bbb K^*\!/\Bbb K^{*3}$ coincides with the image of the
homomorphism $\psi\!:\,E(\Bbb Q)\longrightarrow\tilde E(\Bbb Q)$.
\endproclaim
      The matter is that the kernels of the homomorphisms $\tilde\alpha_{\sssize -}$ and
$\tilde\alpha_{\sssize +}$ do coincide, i\.\,e\. we have a formula similar to the formula
\mythetag{10.22}:
$$
\hskip -2em
\Ker\tilde\alpha_{\sssize -}=\Ker\tilde\alpha_{\sssize +}.
\mytag{10.39}
$$
The formula \mythetag{10.39} is derived from the equality 
$$
\hskip -2em
\tilde\alpha_{\sssize +}(\tilde P)\,\tilde\alpha_{\sssize -}(\tilde P)=1,
\mytag{10.40}
$$
which is similar to \mythetag{10.21}. Indeed, multiplying the formulas \mythetag{10.35} 
and \mythetag{10.37} and applying the curve equation \mythetag{10.2}, we derive  
$$
\hskip -2em
\tilde\alpha_{\sssize +}(\tilde P)\,\tilde\alpha_{\sssize -}(\tilde P)
=\tilde y^2+27\,e^2=\tilde x^3.
\mytag{10.41}
$$
The equality \mythetag{10.41} is equivalent to the equality \mythetag{10.40} since cubic 
factors are neglected modulo $\Bbb K^{*3}$ in $\Bbb K^*\!/\Bbb K^{*3}$.\par
     Let $\tilde P=(\tilde x;\,\tilde y)$ be a non-exceptional rational point of the curve 
\mythetag{10.2} such that $\tilde P\in\Img\psi$. Then its coordinates $\tilde x$ and $\tilde y$ 
are given by the formulas \mythetag{10.4}. Using the second formula \mythetag{10.4} and applying
\mythetag{10.35}, we derive the following relationship:
$$
\hskip -2em
\tilde\alpha_{\sssize -}(\tilde P)=\tilde y-3\,\sqrt{-3\,}\,e
=\frac{y\,x^3-8\,y\,e^2-3\,\sqrt{-3\,}\,e\,x^3}{x^3}.
\mytag{10.42}
$$
Due to the curve equation \mythetag{10.1} the numerator of the fraction in the right 
hand side of \mythetag{10.42} factors as $y\,x^3-8\,y\,e^2-3\,\sqrt{-3\,}\,e\,x^3
=(\tilde y-\sqrt{-3\,}\,e)^3$. Then
$$
\hskip -2em
\tilde\alpha_{\sssize -}(\tilde P)=
\frac{(y-\sqrt{-3\,}\,e)^3}{\tilde x^3}.
\mytag{10.43}
$$
Cubic factors are neglected modulo $\Bbb K^{*3}$ in $\Bbb K^*\!/\Bbb K^{*3}$. 
Therefore the formula \mythetag{10.43} means that $\tilde\alpha_{\sssize -}(\tilde P)=1$ 
and $\tilde P\in\Ker\tilde\alpha_{\sssize -}$.\par
     Thus we have proved the inclusion  $\Img\psi\subseteq\Ker\tilde
\alpha_{\sssize -}$.\par
     Now assume that $\tilde P=(\tilde x;\,\tilde y)$ is some non-exceptional rational 
point of the curve \mythetag{10.2} such that $\tilde P\in\Ker\tilde\alpha_{\sssize -}$. 
Then $\tilde\alpha_{\sssize -}(\tilde P)=1$, which means 
$$
\hskip -2em
\tilde y-3\,\sqrt{-3\,}\,e=z^3\text{, \ where \ }z\neq 0\text{\ \ and \ }z\in\Bbb K.
\mytag{10.44}
$$
The formula \mythetag{10.44} yields $\tilde y=z^3+3\,\sqrt{-3\,}\,e$. Substituting 
$\tilde y=z^3+3\,\sqrt{-3\,}\,e$ into the curve equation \mythetag{10.2}, we derive 
the following equation for $\tilde x$ and $z$: 
$$
\hskip -2em
z^6+6\,\sqrt{-3\,}\,z^3\,e-\tilde x^3=0.
\mytag{10.45}
$$
Using $\tilde x$ and $z$, we define the point $P=(x;\,y)$ of the curve 
\mythetag{10.1} whose coordinates $x$ and $y$ are given by the following 
formulas:
$$
\xalignat 2
&\hskip -2em
x=\frac{2\,\sqrt{-3\,}\,e\,z}{\tilde x-z^2},
&&y=\frac{\sqrt{-3\,}\,e\,(\tilde x+z^2)}{\tilde x-z^2}.
\mytag{10.46}
\endxalignat
$$
The denominator of the fractions \mythetag{10.46} is nonzero. Indeed, otherwise we would
have $\tilde x=z^2$. Substituting $\tilde x=z^2$ into \mythetag{10.45}, we would derive 
$6\,\sqrt{-3\,}\,z^3\,e=0$. Since $e\neq 0$, this would mean $z=0$, which contradicts 
\mythetag{10.44}.\par
     In order to verify that $x$ and $y$ given by the formulas \mythetag{10.46} do 
actually form a point of the curve \mythetag{10.1}, we substitute them into the equation 
\mythetag{10.1}. This yields
$$
\hskip -2em
\frac{4\,e^2\,(z^6+6\,\sqrt{-3\,}\,z^3\,e-\tilde x^3)}{(\tilde x-z^2)^3}=0.
\mytag{10.47}
$$
It is easy to see that the equation \mythetag{10.47} is fulfilled due to \mythetag{10.45}. 
\par     
     The homomorphism $\psi\!:\,E(\Bbb Q)\longrightarrow\tilde E(\Bbb Q)$ is defined 
by means of the formulas \mythetag{10.4}. Substituting \mythetag{10.46} into \mythetag{10.4}, 
we derive 
$$
\gather
\hskip -2em
\tilde x=\frac{6\,\sqrt{-3\,}\,z^3\,e-\tilde x^3+3\,\tilde x^2\,z^2-3\,\tilde x\,z^4+z^6}
{3\,z^2\,(\tilde x-z^2)},
\mytag{10.48}\\
\vspace{1ex}
\hskip -2em
\tilde y=\frac{(3\,\sqrt{-3\,}\,z^3\,e+\tilde x^3-3\,\tilde x^2\,z^2+3\,\tilde x\,z^4-z^6)
\,(\tilde x+z^2)}{3\,z^3\,(\tilde x-z^2)}.
\mytag{10.49}
\endgather
$$
The polynomials in the numerators of the fractions \mythetag{10.48} and \mythetag{10.49} 
are simplified with the use of the equation \mythetag{10.45}. As a result the equality 
\mythetag{10.45} turns to the identity $\tilde x=\tilde x$. The equality \mythetag{10.49} 
turns to the equality $\tilde y=z^3+3\,\sqrt{-3\,}\,e$, which is equivalent to 
\mythetag{10.44}.\par
     The formula \mythetag{10.44} is similar to the formula \mythetag{10.27}. 
Succeeding formulas \mythetag{10.45} through \mythetag{10.49} are similar to the 
formulas \mythetag{10.28} through \mythetag{10.32}. However, there is a crucial difference. 
The number $z$ in \mythetag{10.44} is not rational, while $\tilde x$ and $\tilde y$ are
rational by assumption. There are also explicit irrationalities of the form $\sqrt{-3\,}$ 
in the formulas. For this reason we should take special precautions in order to obtain 
rational numbers $x$ and $y$ in \mythetag{10.46}. Applying \mythetag{10.34} to \mythetag{10.44}, 
we get
$$
\hskip -2em
\tilde y=u^3-9\,u\,v^2+3\,\sqrt{-3\,}\,(u^2\,v-v^3+e).
\mytag{10.50}
$$
Since $\tilde y$, $u$, and $v$ are rational numbers, from \mythetag{10.50} we derive two
equations
$$
\xalignat 2
&\hskip -2em
\tilde y=u^3-9\,u\,v^2,
&&u^2\,v-v^3+e=0. 
\mytag{10.51}
\endxalignat
$$
The first equation \mythetag{10.51} expresses $\tilde y$ through $u$ and $v$. 
\pagebreak The second equation \mythetag{10.51} can be used in order to express 
$e$ through $u$ and $v$:
$$
\hskip -2em
e=v^3-u^2\,v.
\mytag{10.52}
$$
Substituting $\tilde y=u^3-9\,u\,v^2$ and \mythetag{10.52} into the curve equation
\mythetag{10.2}, we derive 
$$
\hskip -2em
u^6+9\,u^4\,v^2+27\,u^2\,v^4+27\,v^6-\tilde x^3=0.
\mytag{10.53}
$$
The equation \mythetag{10.53} replaces the equation \mythetag{10.45}. It factors as follows:
$$
(u^2+3\,v^2-\tilde x)\,(\tilde x^2+\tilde x\,u^2+3\,\tilde x\,v^2+u^4+6\,u^2\,v^2+9\,v^4)=0.
\quad
\mytag{10.54}
$$
Due to \mythetag{10.54} we have two equations which are two options:
$$
\gather
\hskip -2em
u^2+3\,v^2-\tilde x=0,
\mytag{10.55}\\
\tilde x^2+\tilde x\,u^2+3\,\tilde x\,v^2+u^4+6\,u^2\,v^2+9\,v^4=0.
\mytag{10.56}
\endgather
$$
The second option \mythetag{10.56} is a quadratic equation with respect to $\tilde x$. 
One can easily calculate the discriminant of the quadratic equation \mythetag{10.56}:
$$
\hskip -2em
D=-27\,v^4-3\,u^4-18\,u^2\,v^2.
\mytag{10.57}
$$
According to \mythetag{10.44}, the number $z$ is nonzero. Hence its rational 
components $u$ and $v$ in \mythetag{10.34} cannot vanish simultaneously. Then 
\mythetag{10.57} yields the inequality $D<0$. But $\tilde x$ is a rational number by 
assumption, which is incompatible with the inequality $D<0$. Thus we have proved that 
the equality \mythetag{10.55} is the only option for $u$, $v$, and $\tilde x$ derived 
from \mythetag{10.53}. It yields
$$
\hskip -2em
\tilde x=u^2+3\,v^2.
\mytag{10.58}
$$
     Now let's proceed to \mythetag{10.46}. Substituting $z=u+v\,\sqrt{-3\,}$ from 
\mythetag{10.34}, $e=v^3-u^2\,v$ from \mythetag{10.52}, and $\tilde x=u^2+3\,v^2$ from 
\mythetag{10.58} into \mythetag{10.46}, we derive 
$$
\xalignat 2
&\hskip -2em
x=u^2-v^2,
&&y=u^3-u\,v^2,
\mytag{10.59}
\endxalignat
$$
while $\tilde x$ and $\tilde y$ are given by the formulas taken from \mythetag{10.58}
and \mythetag{10.51}:
$$
\xalignat 2
&\hskip -2em
\tilde x=u^2+3\,v^2,
&&\tilde y=u^3-9\,u\,v^2.
\mytag{10.60}
\endxalignat
$$
Thus we have proved that if a point $\tilde P=(\tilde x;\,\tilde y)$ belongs to $\Ker\tilde
\alpha_{\sssize -}$, then there are two rational numbers $u$ and $v$, which are not zero 
simultaneously, such that the coordinates of the point $\tilde P$ are expressed by means 
of the formulas \mythetag{10.60} and there is a point $P=(x;\,y)$, whose coordinates are 
given by the formulas \mythetag{10.59}, such that $\tilde P=\psi(P)$. The latter fact proves 
the inclusion $\Ker\tilde\alpha_{\sssize -}\subseteq\Img\psi$. Combining it with the previously 
proved inclusion $\Img\psi\subseteq\Ker\tilde\alpha_{\sssize -}$, we derive the equality 
$\Img\psi=\Ker\tilde\alpha_{\sssize -}$ which completes the proof of Lemma~\mythelemma{10.12}. 
Lemma~\mythelemma{10.13} is immediate from Lemma~\mythelemma{10.12} due to the equality
\mythetag{10.39}.\par 
\head
11. Factor groups and the rank formula.  
\endhead
     Further steps are similar to those in Section 7. They are visualized by means of 
Fig\.~7.1. The equality \mythetag{7.1} goes without changes:
$$
\hskip -2em
E/\tildepsi(\tilde E(\Bbb Q))\cong(E/\tildepsi\compos\psi(E(\Bbb Q)))\,/
\,(\tildepsi(\tilde E(\Bbb Q))/\tildepsi\compos\psi(E(\Bbb Q))).
\mytag{11.1}
$$
The changes in \mythetag{7.2} are due to the difference of Lemma~\mythelemma{6.4} and
Lemma~\mythelemma{10.5}:
$$
\hskip -2em
[E(\Bbb Q):3\,E(\Bbb Q)]=[E(\Bbb Q):\tildepsi(\tilde E(\Bbb Q))]
\cdot[\tildepsi(\tilde E(\Bbb Q)):\tildepsi\compos\psi(E(\Bbb Q))].
\mytag{11.2}
$$
The formula \mythetag{11.2} follows from \mythetag{11.1} and Lemma~\mythelemma{10.5}.
\par
    Again, we consider the mappings \mythetag{7.3} and derive the isomorphism
\mythetag{7.5} Then, considering the other two isomorphisms \mythetag{7.6} and 
\mythetag{7.7}, we derive the equalities \mythetag{7.8}, \mythetag{7.9}, and
\mythetag{7.10}. Combining them with \mythetag{11.2}, we get
$$
\hskip -2em
[E(\Bbb Q):3\,E(\Bbb Q)]=\frac{[E(\Bbb Q):\tildepsi(\tilde E(\Bbb Q))]
\cdot[\tilde E(\Bbb Q):\psi(E(\Bbb Q))]}
{[\Ker\tildepsi:(\Ker\tildepsi\cap\psi(E(\Bbb Q))]}.
\mytag{11.3}
$$\par
     Unlike the case of Lemma~\mythelemma{7.1}, the denominator of the fraction in
\mythetag{11.3} is always equal to $1$. Indeed, from \mythetag{10.7} we derive that
$\Ker\tildepsi$ is trivial, i\.\,e\. $\Ker\tildepsi=\{\tilde P_\infty\}$. The
intersection $\Ker\tildepsi\cap\psi(E(\Bbb Q))$ is also trivial, i\.\,e\. 
$\Ker\tildepsi\cap\psi(E(\Bbb Q))=\{\tilde P_\infty\}$, which yields 
$[\Ker\tildepsi:(\Ker\tildepsi\cap\psi(E(\Bbb Q))]=1$. As a result \mythetag{11.3}
is written as
$$
\hskip -2em
[E(\Bbb Q):3\,E(\Bbb Q)]=[E(\Bbb Q):\tildepsi(\tilde E(\Bbb Q))]
\cdot[\tilde E(\Bbb Q):\psi(E(\Bbb Q))].
\mytag{11.4}
$$
The next step is to apply Lemmas~\mythelemma{10.8} and \mythelemma{10.12} to \mythetag{11.4}. 
They yield
$$
\hskip -2em
[E(\Bbb Q):3\,E(\Bbb Q)]=|\alpha_{\sssize -}(E(\Bbb Q))|
\cdot|\tilde\alpha_{\sssize -}(\tilde E(\Bbb Q))|.
\mytag{11.5}
$$
Due to \mythetag{10.22} and \mythetag{10.39}, we can write \mythetag{11.5} as
$$
\hskip -2em
[E(\Bbb Q):3\,E(\Bbb Q)]=|\alpha_{\sssize +}(E(\Bbb Q))|
\cdot|\tilde\alpha_{\sssize +}(\tilde E(\Bbb Q))|.
\mytag{11.6}
$$
The formulas \mythetag{11.5} and \mythetag{11.6} are analogs of the formula 
\mythetag{7.19}.\par
     Let's recall that $e$ in \mythetag{10.1} is assumed to be a positive third power free 
integer. Therefore, applying Lemma~\mythelemma{10.1}, we find 
$E_{\sssize\text{tors}}(\Bbb Q)\cong\Bbb Z_3$. The third order elements of $\Bbb Z_3$
do vanish under the tripling endomorphism, i\.\,e\. $3\,E_{\sssize\text{tors}}(\Bbb Q)
=\{\tilde P_\infty\}$ is trivial. Hence we have the following relationships:
$$
\hskip -2em
[E_{\sssize\text{tors}}(\Bbb Q):3\,E_{\sssize\text{tors}}(\Bbb Q)]=|E_3(\Bbb Q)|=
|\Bbb Z_3|=3.
\mytag{11.7}
$$
Applying \mythetag{11.7} to the expression in the left hand sides of \mythetag{11.5} and
\mythetag{11.6}, we get
$$
\hskip -2em
[E(\Bbb Q):3\,E(\Bbb Q)]=3^{\kern 0.5pt r}\!\cdot|E_3(\Bbb Q)|=3^{\kern 0.5pt r+1}\!.
\mytag{11.8}
$$
Combining \mythetag{11.8} with \mythetag{11.5} and \mythetag{11.6}, we can write
these formulas as 
$$
\pagebreak
\hskip -2em
\aligned
&3^{\kern 0.5pt r+1}=|\alpha_{\sssize -}(E(\Bbb Q))|
\cdot|\tilde\alpha_{\sssize -}(\tilde E(\Bbb Q))|,\\
&3^{\kern 0.5pt r+1}=|\alpha_{\sssize +}(E(\Bbb Q))|
\cdot|\tilde\alpha_{\sssize +}(\tilde E(\Bbb Q))|.
\endaligned
\mytag{11.9}
$$
The formulas \mythetag{11.9} reduce the problem of calculating ranks to 
calculating the images of $E(\Bbb Q)$ and $\tilde E(\Bbb Q)$ under the 
descent mappings $\alpha_{\sssize -}$ and $\tilde\alpha_{\sssize -}$ or, 
equivalently, under the descent mappings $\alpha_{\sssize +}$ and 
$\tilde\alpha_{\sssize +}$.\par
\head
12. Images of the descent mappings. 
\endhead
\mylemma{12.1} Let $P$ be a rational point of the curve \mythetag{10.1}. Then 
$\alpha_{\sssize -}(P)$ is presented by some integer number $s=A\cdot B^2$, where 
$A$ and $B$ are two coprime square free positive integer numbers being divisors 
of the number $2\,e$.
\endproclaim
     The case $P=P_\infty$ is trivial. In this case $\alpha_{\sssize -}(P_\infty)=1$
and we choose $A=1$ and $B=1$. If $P=P^{\kern 0.5pt \sssize +}_0$, then the formula
\mythetag{10.9} yields 
$$
\hskip -2em
\alpha_{\sssize -}(P^{\kern 0.5pt \sssize +}_0)=\frac{1}{2\,e}=
\frac{(2\,e)^2}{(2\,e)^3}.
\mytag{12.1}
$$
Cubic factors are neglected modulo $\Bbb Q^{*3}$ in $\Bbb Q^*\!/\Bbb Q^{*3}$. Hence
\mythetag{12.1} is equivalent to $\alpha_{\sssize -}(P^{\kern 0.5pt \sssize +}_0)=
(2\,e)^2$. Let's consider the prime factor expansion of the number $2\,e$:
$$
2\,e=p_1^{\kern 0.4pt\alpha_1}\cdot\ldots\cdot p_r^{\kern 0.4pt\alpha_r}.
\mytag{12.2}
$$
Using $\alpha_1,\,\ldots,\,\alpha_r$ in \mythetag{12.2}, we subdivide 
$p_1,\,\ldots,\,p_r$ into three disjoint sets: 
$$
\align
&\hskip -2em
p_i\in\Cal A\text{\ \ \ if \ \ }2\,\alpha_i\equiv\,1\ \kern -0.6em\mod 3,\\
&\hskip -2em
p_i\in\Cal B\text{\ \ \ if \ \ }2\,\alpha_i\equiv\,2\ \kern -0.6em\mod 3,
\mytag{12.3}\\
&\hskip -2em
p_i\in\Cal C\text{\ \ \ if \ \ }2\,\alpha_i\equiv\,0\ \kern -0.6em\mod 3.
\endalign
$$
Due to \mythetag{12.3} we have the following two coprime square free positive integers:
$$
\xalignat 2
&\hskip -2em
A=\!\!\prod_{p_i\in\Cal A}\!p_i,
&&B=\!\!\prod_{p_i\in\Cal B}\!p_i.
\mytag{12.4}
\endxalignat
$$
Comparing \mythetag{12.4} with \mythetag{12.3} and \mythetag{12.2}, we find that 
the number $(2\,e)^2$ is presented as the product $(2\,e)^2=A\cdot B^2\cdot 
C^{\kern 0.3pt 3}$. Cubic factors are neglected modulo $\Bbb Q^{*3}$ in 
$\Bbb Q^*\!/\Bbb Q^{*3}$. Therefore $\alpha_{\sssize -}(P^{\kern 0.5pt \sssize +}_0)
=A\cdot B^2$ as stated in Lemma~\mythelemma{12.1}.\par
     If the point $P$ is different from $P_\infty$ and $P^{\kern 0.5pt \sssize +}_0$, 
we consider the coordinates of this point: $P=(x;\,y)$. Since $P$ is a rational point,
we have two irreducible fractions
$$
\xalignat 2
&\hskip -2em
x=\frac{m}{q_1},
&&y=\frac{n}{q_2}.
\mytag{12.5}
\endxalignat
$$
Substituting the fractions \mythetag{12.5} into the curve equation \mythetag{10.1}, we 
derive
$$
\hskip -2em
\frac{m^3}{q_1^3}=\frac{n^2-e^2\,q_2^2}{q_2^2}.
\mytag{12.6}
$$
Both sides of the equality \mythetag{12.6} are  irreducible fractions. Hence $q_1^3=q_2^2$.
From the equality $q_1^3=q_2^2$ we derive that the numbers $q_1$ and $q_2$ are presented as 
$$
\xalignat 2
&\hskip -2em
q_1=q^2,
&q_2=q^3
\mytag{12.7}
\endxalignat
$$
(compare with \mythetag{8.5}). Applying \mythetag{12.7} to \mythetag{12.6}, we obtain the
equality 
$$
\hskip -2em
m^3=(n-e\,q^3)\,(n+e\,q^3), 
\mytag{12.8}
$$
while the formulas \mythetag{12.5} for the coordinates $x$ and $y$ turn to the following ones:
$$
\xalignat 2
&\hskip -2em
x=\frac{m}{q^2},
&&y=\frac{n}{q^3}.
\mytag{12.9}
\endxalignat
$$\par
     Now let's recall, the formula \mythetag{10.8}. Using \mythetag{12.9}, this formula yields
$$
\hskip -2em
\alpha_{\sssize -}(P)=y-e=\frac{n-e\,q^3}{q^3}. 
\mytag{12.10}
$$
Cubic factors are neglected modulo $\Bbb Q^{*3}$ in $\Bbb Q^*\!/\Bbb Q^{*3}$. Hence
we can write \mythetag{12.10} as 
$$
\hskip -2em
\alpha_{\sssize -}(P)=s=n-e\,q^3. 
\mytag{12.11}
$$
Comparing \mythetag{12.11} with \mythetag{12.8}, we obtain the following Diophantine
equation 
$$
\hskip -2em
m^3=s\,(s+2\,e\,q^3). 
\mytag{12.12}
$$
The fractions \mythetag{12.9} are irreducible. Therefore $\gcd(n,q)=1$. Applying this
equality to \mythetag{12.11}, we derive $\gcd(s,q)=1$. Keeping in mind $\gcd(s,q)=1$, 
let's calculate the greatest common divisor of two multiplicands $s$ and 
$s+2\,e\,q^3$ in the right hand side of the equality \mythetag{12.12} and denote it 
through $M$: 
$$
\hskip -2em
M=\gcd(s,s+2\,e\,q^3)=\gcd(s,2\,e\,q^3)=\gcd(s,2\,e). 
\mytag{12.13}
$$
From \mythetag{12.13} it is clear that $M$ is a divisor of the number $2\,e$.
\par
     Like in \mythetag{12.2}, let's consider the prime factor expansion of the number
$m$ from the left hand side of the equation \mythetag{12.12}:
$$
m=\pm\,p_1^{\kern 0.4pt\gamma_1}\cdot\ldots\cdot p_\rho^{\kern 0.4pt\gamma_\rho}.
\mytag{12.14}
$$
If $p_i$ in \mythetag{12.14} is not a divisor of $M$ in \mythetag{12.13}, then due to
\mythetag{12.12} either $p_i^{\,3\,\gamma_i}$ divides $s$ or $p_i^{\,3\,\gamma_i}$ 
divides $s+2\,e\,q^3$. The case where $p_i$ is a divisor of $M$ is more complicated. 
In this case $3\,\gamma_i=\alpha_i+\beta_i$, where $\alpha_i\neq 0$, $\beta_i\neq 0$, 
so that $p_i^{\,\alpha_i}$ divides $s$, and $p_i^{\,\beta_i}$ divides $s+2\,e\,q^3$. 
Prime factors of this sort are grouped into three subsets: 
$$
\align
&\hskip -2em
p_i\in\Cal A\text{\ \ \ if \ \ }\alpha_i\equiv\,1\ \kern -0.6em\mod 3
\text{\ \ \ and \ \ }\beta_i\equiv\,2\ \kern -0.6em\mod 3,\\
&\hskip -2em
p_i\in\Cal B\text{\ \ \ if \ \ }\alpha_i\equiv\,2\ \kern -0.6em\mod 3
\text{\ \ \ and \ \ }\beta_i\equiv\,1\ \kern -0.6em\mod 3,
\mytag{12.15}\\
&\hskip -2em
p_i\in\Cal C\text{\ \ \ if \ \ }\alpha_i\equiv\,0\ \kern -0.6em\mod 3
\text{\ \ \ and \ \ }\beta_i\equiv\,0\ \kern -0.6em\mod 3.
\endalign
$$
Using \mythetag{12.15}, we define two coprime square free positive 
integer numbers:
$$
\xalignat 2
&\hskip -2em
A=\!\!\prod_{p_i\in\Cal A}\!p_i,
&&B=\!\!\prod_{p_i\in\Cal B}\!p_i.
\mytag{12.16}
\endxalignat
$$
From \mythetag{12.15} and \mythetag{12.16} we immediately derive \pagebreak 
$s=A\cdot B^2\cdot C^{\kern 0.3pt 3}$. Cubic factors are neglected modulo 
$\Bbb Q^{*3}$ in $\Bbb Q^*\!/\Bbb Q^{*3}$. Therefore we can write 
$\alpha_{\sssize -}(P)=s=A\cdot B^2$. The numbers $A$ and $B$ are divisors 
of $M$ by construction, while $M$ is a divisor of $2\,e$. Lemma~\mythelemma{12.1}
is proved.\par
       Lemma~\mythelemma{12.1} yields a very rough description of the set 
$\alpha_{\sssize -}(E(\Bbb Q))$. It can be refined in the following way. 
\mylemma{12.2} An integer number $s=A\cdot B^2$, where $A$ and $B$ are two coprime square 
free positive integers dividing $2\,e$, represents an element of\/ 
$\alpha_{\sssize -}(E(\Bbb Q))$ if and only if the following Diophantine equation 
has a non-trivial solution $X$, $Y$, $Z$:
$$
\hskip -2em
(A\cdot B^2)\,X^3+(A^2\cdot B)\,Y^3+(2\,e)\,Z^3=0.
\mytag{12.17}
$$
\endproclaim
\demo{Proof} Necessity. If $s=\alpha_{\sssize -}(P_\infty)=1$, then $A=1$ and $B=1$.
In this case the equation \mythetag{12.17} is solvable and $X=1$, $Y=-1$, $Z=0$ is 
its explicit solution.\par 
      If $s=\alpha_{\sssize -}(P^{\kern 0.5pt \sssize +}_0)$, then, as we have seen in
proving Lemma~\mythelemma{12.1}, $s=A\cdot B^2$, while $A\cdot B^2\cdot C^{\kern 0.3pt 3}
=(2\,e)^2$. Multiplying both sides of the equality $A\cdot B^2\cdot C^{\kern 0.3pt 3}
=(2\,e)^2$ by $A^2\cdot B\cdot(2\,e)$, we derive the following equality:
$$
\hskip -2em
(A^2\cdot B)\,(2\,e)^3=A^2\cdot B\cdot(2\,e)\cdot A\cdot B^2\cdot C^{\kern 0.3pt 3}=
(2\,e)\,(A\,B\,C)^3.
\mytag{12.18}
$$
The equality \mythetag{12.18} yields the solution of the Diophantine equation \mythetag{12.17} 
with $X=0$, $Y=2\,e\neq 0$, $Z=-A\,B\,C$.\par
     Now assume that $s=\alpha_{\sssize -}(P)$ for a rational point $P$ of the curve
\mythetag{10.1} other than $P_\infty$ and $P^{\kern 0.5pt \sssize +}_0$. As we have seen
above in proving Lemma~\mythelemma{12.1}, in this case $s$ is a solution of the equation
\mythetag{12.12} given by the formula 
$$
\hskip -2em
s=A\cdot B^2\cdot C^{\kern 0.3pt 3}.
\mytag{12.19}
$$ 
This formula is derived from \mythetag{12.15} and \mythetag{12.16}. Similarly, from 
\mythetag{12.15} and \mythetag{12.16} one can derive the following formula for
$s+2\,e\,q^3$:
$$
\hskip -2em
s+2\,e\,q^3=A^2\cdot B\cdot D^3. 
\mytag{12.20}
$$
Subtracting \mythetag{12.19} from \mythetag{12.20}, we obtain the equality
$$
\hskip -2em
(A^2\cdot B)\,D^3-(A\cdot B^2)\,C^{\kern 0.3pt 3}=(2\,e)\,q^3. 
\mytag{12.21}
$$
It is easy to see that the equality \mythetag{12.21} yields the solution of the 
Diophantine equation \mythetag{12.17} with $X=-C$, $Y=D$, $Z=-q\neq 0$.\par
     Sufficiency. Assume that $X$, $Y$, $Z$ are three integer numbers composing 
a non-trivial solution of the equation \mythetag{12.17}. If $Z=0$, then $X\neq 0$
or $Y\neq 0$. In this case \mythetag{12.17} reduces to $B\,X^3+A\,Y^3=0$. It yields
the implications
$$
\hskip -2em
\gathered
X\neq 0\Longrightarrow Y^3=-\frac{B}{A}\,X^3\Longrightarrow Y\neq 0,\\
\vspace{1ex}
Y\neq 0\Longrightarrow X^3=-\frac{A}{B}\,X^3\Longrightarrow X\neq 0.
\endgathered
\mytag{12.22}
$$
Using the equation $B\,X^3+A\,Y^3=0$ and taking into account \pagebreak
the implications \mythetag{12.22}, we can transform the formula $s=A\cdot B^2$ 
in the following way:
$$
\hskip -2em
s=A\cdot B^2=-\frac{Y^3}{X^3}\,B^3=\biggl(-\frac{Y\,B}{X}\biggr)^{\lower 2pt
\hbox{$\kern -0.5pt\ssize 3$}}\!. 
\mytag{12.23}
$$
Cubic factors are neglected modulo $\Bbb Q^{*3}$ in $\Bbb Q^*\!/\Bbb Q^{*3}$. The
formula \mythetag{12.23} means that $s=A\cdot B^2$ in this case is equivalent to
$s=1$, which is the value of the mapping $\alpha_{\sssize -}$ at the point 
$P=P_\infty$.\par
     Now assume that $Z\neq 0$. If under this assumption $Y=0$, then the equation 
\mythetag{12.17} yields $X\neq 0$ and provides the following equalities:
$$
\hskip -2em
s=A\cdot B^2=-(2\,e)\,\frac{Z^3}{X^3}. 
\mytag{12.24}
$$ 
Cubic factors are neglected modulo $\Bbb Q^{*3}$ in $\Bbb Q^*\!/\Bbb Q^{*3}$. The
formula \mythetag{12.24} means that $s=A\cdot B^2$ in this case is equivalent to
$s=-2\,e$, which is the value of the mapping $\alpha_{\sssize -}$ at the point 
$P=P^{\kern 0.5pt \sssize -}_0$.\par
     The case $Z\neq 0$ and $X=0$ is another option. In this case the equation 
\mythetag{12.17} yields $Y\neq 0$ and provides the following equalities:
$$
\hskip -2em
s=A\cdot B^2=\frac{(A^2\cdot B)^2}{A^3}=\frac{1}{A^3}
\biggl(-2\,e\,\frac{Z^3}{Y^3}\biggr)^{\lower 2pt
\hbox{$\kern -0.5pt\ssize 2$}}=\frac{1}{2\,e}
\biggl(\frac{2\,e\,Z^2}{A\,Y^2}\biggr)^{\lower 2pt
\hbox{$\kern -0.5pt\ssize 3$}}\!.
\mytag{12.25}
$$ 
Cubic factors are neglected modulo $\Bbb Q^{*3}$ in $\Bbb Q^*\!/\Bbb Q^{*3}$. 
The formula \mythetag{12.25} means that $s=A\cdot B^2$ in this case is equivalent to
$s=1/(2\,e)$, which is the value of the mapping $\alpha_{\sssize -}$ at the point 
$P=P^{\kern 0.5pt \sssize +}_0$.\par
      In the general case we have a solution of the equation \mythetag{12.17} with
$X\neq 0$, $Y\neq 0$, $Z\neq 0$. In this case we can write the equation \mythetag{12.17}
as follows:
$$
\hskip -2em
(A^2\cdot B)\,Y^3=(A\cdot B^2)\,(-X)^3+(2\,e)\,(-Z)^3.
\mytag{12.26}
$$
Cubic factors are neglected modulo $\Bbb Q^{*3}$ in $\Bbb Q^*\!/\Bbb Q^{*3}$. Therefore
$s=A\cdot B^2$ is equivalent to $s=A\cdot B^2\cdot (-X)^3$. For the sake of similarity
to \mythetag{12.12} let's denote $q=-Z$. As a result the above equality \mythetag{12.26} 
turns to
$$
\hskip -2em
(A^2\cdot B)\,Y^3=s+2\,e\,q^3.
\mytag{12.27}
$$
Multiplying both sides of \mythetag{12.27} by $s=A\cdot B^2\cdot (-X)^3$, we find that 
the product $s\,(s+2\,e\,q^3)$ is an exact cube. Indeed, we have
$$
\hskip -2em
s\,(s+2\,e\,q^3)=(-A\,B\,X\,Y)^3. 
\mytag{12.28}
$$
Denoting $m=-A\,B\,X\,Y$ in \mythetag{12.28}, we obtain the equality coinciding with
\mythetag{12.12}.
$$
\hskip -2em
s\,(s+2\,e\,q^3)=m^3. 
\mytag{12.29}
$$
The rest is to denote $n=s+e\,q^3$ and compose two fractions:
$$
\xalignat 2
&\hskip -2em
x=\frac{m}{q^2},
&&y=\frac{n}{q^3}.
\mytag{12.30}
\endxalignat
$$
Due to \mythetag{12.29} the number $n=s+e\,q^3$ satisfies the equation \mythetag{12.8},
while the values of the fractions \mythetag{12.30} satisfy the curve equation 
\mythetag{10.1}. They define a rational point $P=(x;\,y)$ of this curve different from
$P_\infty$, $P^{\kern 0.5pt \sssize +}_0$, and $P^{\kern 0.5pt \sssize -}_0$. Applying
he formula \mythetag{10.8} to this point, we derive the formula for $\alpha_{\sssize -}(P)$:
$$
\hskip -2em
\alpha_{\sssize -}(P)=y-e=\frac{n-e\,q^3}{q^3}=
\frac{A\cdot B^2\cdot (-X)^3}{(-Z)^3}.
\mytag{12.31}
$$
Cubic factors are neglected modulo $\Bbb Q^{*3}$ in $\Bbb Q^*\!/\Bbb Q^{*3}$. Therefore
the formula \mythetag{12.31} means that $s=A\cdot B^2$ is a value of the mapping 
$\alpha_{\sssize -}$ in the general case too. Lemma~\mythelemma{12.2} is proved.
\qed\enddemo
     The mapping $\tilde\alpha_{\sssize -}$ is somewhat different from the mapping 
$\alpha_{\sssize -}$. Its domain is the set of rational points of the second curve 
\mythetag{10.2}. In defining this mapping we replaced the field of rational numbers 
$\Bbb Q$ by its algebraic extension $\Bbb K=\Bbb Q(\sqrt{-3\,})$. The following 
definitions are standard.
\mydefinition{12.1} A polynomial of one variable is called monic if its leading
coefficient is equal to unity: $f(x)=x^n+a_1\,x^{n-1}+\ldots+a_{n-1}\,x+a_n$.
\enddefinition
\mydefinition{12.2} A number $x$ from some extension of the field of rational numbers 
$\Bbb Q$ is called algebraic if it is a root of some polynomial with rational 
coefficients.
\enddefinition
\mydefinition{12.3} An algebraic number $x$ is called an algebraic integer 
if it is a root of some monic polynomial with integer coefficients.
\enddefinition
\mydefinition{12.4} An algebraic integer number $x$ belonging to the number field 
$\Bbb K=\Bbb Q(\sqrt{-3\,})$ is called an Eisenstein integer (see \mycite{16}).
\enddefinition
     Eisenstein integers constitute a ring within the field $\Bbb K=\Bbb Q(\sqrt{-3\,})$.
We denote this ring through $\Bbb Z(\sqrt{-3\,})$. The ring of Eisenstein integers
$\Bbb Z(\sqrt{-3\,})$ is an integral domain, i\.\,e\. it has no divisors of zero. The
field $\Bbb K=\Bbb Q(\sqrt{-3\,})$ coincides with the field of fractions for the 
ring $\Bbb Z(\sqrt{-3\,})$. 
\mydefinition{12.5} An algebraic integer number $x\in\Bbb Z(\sqrt{-3\,})$ is called
invertible if its inverse element $x^{-1}=1/x$ is also an algebraic integer, i\.\,e\.
if $x^{-1}\in\Bbb Z(\sqrt{-3\,})$.
\enddefinition
\noindent
There are exactly six invertible elements in the ring of Eisenstein integers
$\Bbb Z(\sqrt{-3\,})$:
$$
\xalignat 3
&\hskip -2em
\varepsilon=\frac{1+\sqrt{-3\,}}{2},
&&\omega=\frac{-1+\sqrt{-3\,}}{2},
&&-1,\quad\\
\vspace{-1.5ex}
\mytag{12.32}
\\
\vspace{-1.5ex}
&\hskip -2em
\overline\omega=\frac{-1-\sqrt{-3\,}}{2},
&&\overline\varepsilon=\frac{1-\sqrt{-3\,}}{2},
&&\hphantom{-}\,\,1.\quad
\endxalignat
$$
Here $\varepsilon$ is an elementary sixth root of unity, while $\omega$ is an
elementary cubic root of unity. The ring of regular integers
$\Bbb Z$ has only two invertible elements: $1$ and $-1$.\par
      Eisenstein integers from the ring $\Bbb Z(\sqrt{-3\,})$ constitute a 
two-dimensional grid. Indeed, each element $x\in\Bbb Z(\sqrt{-3\,})$ is presented as 
$$
\hskip -2em
x=u+v\,\varepsilon\text{, \ where \ }u\in\Bbb Z,\ v\in\Bbb Z,
\mytag{12.33}
$$
and where $\varepsilon$ is taken from \mythetag{12.32}. The formula \mythetag{12.33}
is similar to \mythetag{10.34}.\par
\mydefinition{12.6} An algebraic integer number $x$ from the ring of Eisenstein integers 
$\Bbb Z(\sqrt{-3\,})$ is called prime if it is non-invertible and if it cannot be presented 
as a product of two other non-invertible algebraic integers from this ring. 
\enddefinition
     It is very important that the ring of Eisenstein integers $\Bbb Z(\sqrt{-3\,})$
is a Euclidean domain (see \mycite{16}), where the norm of the number \mythetag{12.33} 
is given by the formula
$$
|u+v\,\varepsilon|^2=u^2+u\,v\,+v^2.
$$ 
Each Euclidean domain is a unique factorization domain. When applied to the ring of
Eisenstein integers, this means that each Eisenstein integer $x\in\Bbb Z(\sqrt{-3\,})$ 
has an expansion into the product of Eisenstein primes:
$$
\hskip -2em
x=p_1\cdot\ldots\cdot p_r
\mytag{12.34}
$$
The expansion \mythetag{12.34} is unique up to the order of multiplicands and 
up to multiplying the prime numbers $p_1,\,\ldots,\,p_r$ by invertible elements 
\mythetag{12.32}. Due to expansions of the form \mythetag{12.34} one can 
formulate two lemmas similar to Lemma~\mythelemma{12.1} and Lemma~\mythelemma{12.2}.
\mylemma{12.3} Let $\tilde P$ be a rational point of the curve \mythetag{10.2}. Then 
$\tilde\alpha_{\sssize -}(\tilde P)$ is presented by some Eisenstein integer 
$s=\eta\,\cdot A\cdot B^2$, where $A$ and $B$ are two nonzero coprime square free Eisenstein 
integers being divisors of the number $6\,\sqrt{-3\,}\,e$ and where $\eta$ is one of 
the three invertible Eisenstein integers $1$, $\varepsilon$, or $\omega$ from 
\mythetag{12.32}.
\endproclaim
\mylemma{12.4} An Eisenstein integer number $s=\eta\,\cdot A\cdot B^2$, where $A$ and 
$B$ are two nonzero coprime square free Eisenstein integers being divisors of 
$6\,\sqrt{-3\,}\,e$ and where $\eta$ is one of the three invertible Eisenstein integers 
$1$, $\varepsilon$, or $\omega$ from \mythetag{12.32}, represents an element of 
$\tilde\alpha_{\sssize -}(\tilde E(\Bbb Q))$ if and only if the homogeneous cubic 
equation 
$$
\hskip -2em
(\eta\,\cdot A\cdot B^2)\,X^3+(\eta^{-1}\!\cdot A^2\cdot B)\,Y^3
+(6\,\sqrt{-3\,}\,e)\,Z^3=0
\mytag{12.35}
$$
has a non-trivial solution $X$, $Y$, $Z$ in $\Bbb Z(\sqrt{-3\,})$ such that if 
$Z\neq 0$, then the values of the following two fractions are 
regular rational numbers:
$$
\xalignat 2
&x=\frac{-A\,B\,X\,Y}{Z^2},
&&y=\frac{(\eta\,\cdot A\cdot B^2)\,X^3+(3\,\sqrt{-3\,}\,e)\,Z^3}{Z^3}.
\endxalignat
$$
\endproclaim
Note that the number $6\,\sqrt{-3\,}$ used in \mythetag{12.35} and in the statements of
the above lemmas belongs to the ring of Eisenstein integers $\Bbb Z(\sqrt{-3\,})$. Indeed, 
we have the following presentation of the form \mythetag{12.33} for this number:
$$
6\,\sqrt{-3\,}=-6+12\,\varepsilon.
$$
As for the proofs of Lemma~\mythelemma{12.3} and Lemma~\mythelemma{12.4}, they use 
almost the same arguments as the proofs of Lemma~\mythelemma{12.1} and 
Lemma~\mythelemma{12.2}.\par
\head
13. Conclusions.
\endhead
     Elliptic curves \mythetag{1.1} brought to the form \mythetag{2.2} are
associated with perfect cuboids. Potentially they could be used in finding an
example of such a cuboid or in proving their non-existence. \pagebreak 
Therefore our further efforts will be directed to describing rational points 
of these curves by some formula\kern 0.4em or, which is more likely, to 
evaluating them numerically. Computing the ranks of the curves is the first 
step in this direction. For the case where $4\,R^{\kern 0.5pt 2}\,N$ in 
\mythetag{2.2} is an exact cube the rank problem is already solved by 
Theorem~\mythetheorem{5.7} within the $2$-descent method.\par
     The case where $4\,R^{\kern 0.5pt 2}\,N$ is not an exact cube is more 
complicated. In this case the formulas \mythetag{11.9} and Lemmas~\mythelemma{12.1}, 
\mythelemma{12.2}, \mythelemma{12.3}, \mythelemma{12.4} obtained within the 
$3$-descent method could be a background for some numerical algorithm. Building 
such an algorithm and analyzing its output is a subject for a separate article. 
\head
14. Acknowledgments.
\endhead
     We are grateful to \myhref{http://e-science.ru/forum/index.php?showuser=28577}
{Sonic86} who gave us the reference to \mycite{7} on 
\myhref{http://e-science.ru/forum/index.php?showtopic=40240}{e-science.ru} forum.

\enddocument
\end